\documentclass{svmult}
\usepackage{bbm}
\usepackage{amsfonts}
\usepackage{mathrsfs}
\usepackage{amssymb}
\usepackage{cases}
\newenvironment{enumerate-roman}{\begin{enumerate}}{\end{enumerate}}
\usepackage{makeidx}         
\usepackage{graphicx}        
\usepackage{multicol}        
\usepackage[bottom]{footmisc}

\makeindex             

\newtheorem{thm}{Theorem}[section]
\newtheorem{lem}[thm]{Lemma}
\newtheorem{prop}[thm]{Proposition}
\newtheorem{cor}[thm]{Corollary}
\newtheorem{rmk}[thm]{Remark}

\newtheorem{defi}[thm]{Definition}

\begin{document}

\baselineskip14pt

\title*{Probabilistic Representation of Weak Solutions of Partial Differential Equations with Polynomial Growth Coefficients\thanks{Accepted by Journal of Theoretical Probability, to appear. DOI: 10.1007/s10959-011-0350-y}}
\titlerunning{Weak Solutions of PDEs with p-Growth Coefficients}
\author{Qi Zhang\inst{\ {\rm a}}, Huaizhong Zhao\inst{\ {\rm b}}}
\authorrunning{Q. Zhang and H.Z. Zhao}
\institute{$^{\rm a}$ School of Mathematical Sciences, Fudan
University, Shanghai, 200433, China.\\
$^{\rm b}$ Department of Mathematical Sciences,
Loughborough University, Loughborough, LE11 3TU, UK.\\
\texttt{Emails: qzh@fudan.edu.cn}; \texttt{H.Zhao@lboro.ac.uk}}

\maketitle
\newcounter{bean}
\begin{abstract}
In this paper we develop a new weak convergence and compact
embedding method to study the existence and uniqueness of the
$L_{\rho}^{2p}({\mathbb{R}^{d}};{\mathbb{R}^{1}})\times
L_{\rho}^2({\mathbb{R}^{d}};{\mathbb{R}^{d}})$ valued solution of
backward stochastic differential equations with p-growth
coefficients. Then we establish the probabilistic representation of
the weak solution of PDEs with p-growth coefficients via
corresponding BSDEs.
\end{abstract}

\vskip5pt

\vskip5pt

\renewcommand{\theequation}{\arabic{section}.\arabic{equation}}

\section{Introduction}
In this paper, we study the probabilistic representation of the weak
solution of a class of parabolic partial differential equations
(PDEs) on $\mathbb{R}^{d}$ with p-growth coefficients
\begin{numcases}{}\label{1jan2009}
{{\partial v}\over{\partial
t}}(t,x)=\mathscr{L}v(t,x)+f\big(x,v(t,x),(\sigma^*\nabla
v)(t,x)\big),\ \ \ \ 0\leq t\leq
T,\nonumber\\
v(0,x)=h(x),
\end{numcases}
by the solution of the corresponding backward stochastic
differential equations (BSDEs) in $\rho$-weighted $L^2$ space. Here
$\mathscr{L}$ is a second order differential operator
\begin{eqnarray}\label{h}
\mathscr{L}={1\over2}\sum_{i,j=1}^da_{ij}(x){{\partial^2}\over{\partial
x_i\partial x_j}}+\sum_{i=1}^db_i(x){\partial\over{\partial x_i}},
\end{eqnarray}
$(a_{ij}(x))$ is a symmetric matrix with a decomposition
$(a_{ij}(x))=(\sigma_{ij}(x))(\sigma_{ij}(x))^*$, $f: (x,y,z)\mapsto
f(x,y,z)$ is a function of polynomial growth in $y$ and Lipschitz
continuous in $z$. Many partial differential equations arising in
physics, engineering and biology have polynomial growth nonlinear
terms e.g. KPP-Fisher equations, Allen-Cahn equations and
Ginzburg-Landau equations.
The representation provides an important
connection between stochastic flows generated by $\mathscr{L}$ and
the weak solutions of PDEs possibly with polynomial growth
coefficients. In connection with the classical solutions of the
linear parabolic PDEs, the well-known Feynman-Kac formula provides
the probabilistic representation for them and originated many
important developments (Feynman \cite{fe}, Kac \cite{ka}). An
alternative probabilistic representation using only the values of a
finite (random) set of times to the linear heat equations was
obtained recently by Dalang, Mueller and Tribe \cite{dmt}. This idea
made it possible for them to obtain corresponding formula for a
wide class of linear PDEs such as some wave equations with
potentials. The Feynman-Kac formula has played important roles in
problems such as the large deviation theory of Donsker and Varadhan
\cite{do-va}, Wentzell and Freidlin \cite{we-fr}, small time
asymptotics of heat kernel and its logarithmic derivatives, in
particular on Riemannian manifolds (Elworthy \cite{el}, Malliavin
and Stroock \cite{ma-st}). The Feynman-Kac formula has been extended
and used to quasi-linear parabolic type partial differential
equations, especially, in the study of the generalized KPP equations
using the large deviation theory method by Freidlin \cite{fr}, using
the semi-classical probabilistic method by Elworthy, Truman and Zhao
\cite{el-tr-zh-ga}.  The study of the quasi-linear parabolic type
PDE is based on an equation of the Feynman-Kac type integration of
stochastic functionals. The approach of the backward stochastic
differential equations, pioneered by Pardoux and Peng \cite{pa-pe1},
\cite{pa-pe2} originally, provided an alternative approach to the
classical solution of the parabolic type PDEs, when the coefficients
of the PDE are sufficiently regular and Lipschitz continuous. This
was extended to the viscosity solution of a large class of partial
differential equations and BSDEs. They include the linear growth
case considered by Lepeltier and San Martin \cite{le-sa}, the
quadratic coefficients (in $z$) considered by Kobylanski \cite{ko},
Briand and Hu \cite{br-hu}, and the polynomial growth coefficients
in Pardoux \cite{pa1}. The solution of the BSDEs in above cases
gives the probabilistic representation of the classical or viscosity
solution of the PDEs as a generalization to the Feynman-Kac formula.
Applications of BSDEs have been found in some problems such as a
model in mathematics of finance (El Karoui, Peng and Quenez
\cite{el-pe-qu}), as an efficient method for constructing
$\Gamma$-martingales on Riemannian manifolds (Darling \cite{da}),
and as an intrinsic tool to construct the pathwise stationary
solution for stochastic PDEs (Zhang and Zhao \cite{zh-zh1},
\cite{zh-zh2}).

The Feynman-Kac approach to a Sobolev or $L^2$ space valued weak
solution of PDEs has been concentrated mainly on linear problems.
Many important progress has been made e.g. in quantum field theory
(see \cite{si}). The probabilistic approach to the weak solution of
quasi-linear PDEs stayed  behind. Regularity of the solutions, even
in the sense of weak derivative, was not given in Freidlin's
probabilistic approach of generalized solution formally represented
by the Feynman-Kac formula (\cite{fr}).  The BSDEs start to show
some usefulness in this aspect, when
 the coefficients are of Lipschitz continuous in the
space $L_{\rho}^2({\mathbb{R}^{d}};{\mathbb{R}^{1}})\times
L_{\rho}^2({\mathbb{R}^{d}};{\mathbb{R}^{d}})$ or of linear growth,
and monotone,  from the work of Barles and Lesigne \cite{ba-le},
Bally and Matoussi \cite{ba-ma}, Zhang and Zhao \cite{zh-zh1},
\cite{zh-zh2}. The objective of this paper is to move away from the
assumption of the linear growth of $f$ and from considering the
classical or viscosity solution of PDEs to establish the
probabilistic representation for the weak solution of such
polynomial growth PDEs. Although the connection of BSDEs with the
viscosity solution for the cases of quadratic and polynomial growth
has been obtained in \cite{ko}, \cite{pa1} respectively, the
existing methods in the study of BSDEs for finding the solution of
the BSDEs in $L_{\rho}^2({\mathbb{R}^{d}};{\mathbb{R}^{1}})\times
L_{\rho}^2({\mathbb{R}^{d}};{\mathbb{R}^{d}})$ are not adequate to
solve the problem of the weak solution of BSDEs with p-growth
coefficients. The fixed point method in
$M^{2}([t,T];L_{\rho}^2({\mathbb{R}^{d}};{\mathbb{R}^{1}}))\times
M^{2}([t,T];L_{\rho}^2({\mathbb{R}^{d}};{\mathbb{R}^{d}}))$, which
is equivalent to finding a strongly convergent sequence in the same
space, seems difficult to work for the problem with p-growth
coefficients. It is also inadequate to use a combination of the weak
convergence in finite dimensional space developed by Pardoux
\cite{pa1} and the weak solution method developed by Bally and
Matoussi \cite{ba-ma}, Zhang and Zhao \cite{zh-zh1}, \cite{zh-zh2},
to solve this problem. We need to introduce some new ideas to the
study of BSDEs. The progress of this problem was made when we
realized that, in addition to 
the method of Zhang and Zhao (\cite{zh-zh1}, \cite{zh-zh2}), as well
as the standard approach using Alaoglu lemma to find a weakly
convergent sequence $(Y^n,Z^n)$, we can use the equivalence of norm
principle and Rellich-Kondrachov Compactness Theorem to get a
strongly convergent sequence $Y^n$.
Our recent result on the $S^2([t,T],
L_{\rho}^2({\mathbb{R}^{d}};{\mathbb{R}^{1}}))\times
M^{2}([t,T];L_{\rho}^2({\mathbb{R}^{d}};{\mathbb{R}^{d}}))$ valued
solution of BDSDEs with non-Lipschitz linear growth coefficients made
it possible for us to study the BSDEs in $S^{2p}([t,T],
L_{\rho}^{2p}({\mathbb{R}^{d}};{\mathbb{R}^{1}}))\times
M^{2}([t,T];L_{\rho}^2({\mathbb{R}^{d}};{\mathbb{R}^{d}}))$ with
polynomial growth coefficients, even without assuming $f$ being
locally Lipschitz continuous in $y$. Of course, we need to assume
the monotonicity condition of $f$ in $y$. Moreover, it is also an
essential step to prove the strong convergence of $Z^n$ in
$M^2([t,T], L_{\rho}^2({\mathbb{R}^{d}};{\mathbb{R}^{d}}))$ from the
result of the strong convergence of $Y^n$ and It$\hat {\rm o}$'s
formula. The weak convergence and compact embedding method has been
used in the study of PDEs. However, as far as we know, to use this
kind of argument to the study of BSDE, this paper is the first time
in literature. The equivalence of norm principle and very careful
probabilistic (measure theoretical) and analytic arguments including
localization made it work in the probabilistic context. However, the
probabilistic case is a lot more complicated than the deterministic
PDEs case as we need to work on the space
$\Omega\times[0,T]\times\mathbb{R}^d$ and solve the equation with
probability one, instead only work on $[0,T]\times\mathbb{R}^d$
in the deterministic PDEs case. The probabilistic representation can
be regarded as a generalized Feynman-Kac formula to the weak
solution of the PDEs with p-growth coefficients. We would like to point out that analysts already studied PDE with polynomial growth coefficients when these coefficients do not depend on $\nabla u$ (see Robinson \cite{ro}, Temam \cite{te}). Our method pushes the study to more general equation with nonlinear term depending on $\nabla u$. But this is not the main purpose of the paper. Our main purpose is to find a method to solve the BSDEs in $L_{\rho}^{2p}({\mathbb{R}^{d}};{\mathbb{R}^{1}})\times L_{\rho}^2({\mathbb{R}^{d}};{\mathbb{R}^{d}})$ space with polynomial growth coefficients, so it provides a probabilistic representation to the corresponding PDEs. This is new in literature. Moreover, our approach does not depend on results of PDEs. Rather we can obtain results about PDEs from the study of BSDEs. Due to this important aspect of our results here, we can extend this result to backward doubly stochastic differential equations (BDSDEs) so that we obtain new results on the stationary solutions of SPDEs via BDSDEs with polynomial growth coefficients. See Zhang and Zhao \cite{zh-zh4}. We believe our method will be useful to other types of PDEs or SPDEs and BSDEs or BDSDEs as well.

After this paper was completed, we were informed the paper Matoussi
and Xu \cite{MX}. But we would like to point out what we have proved
as well as our methods are different. Notice the convergence
$(Y_s^{t,x,n}, Z_{s}^{t,x,n})$ is only a weak convergence along a
subsequence according to the Alaoglu lemma. If one considers weak
convergence in $M^2([t,T], {\mathbb{R}^{1}})\times
M^{2}([t,T];{\mathbb{R}^{d}}))$, which worked well in Pardoux
\cite{pa1} for the case of viscosity solutions of the PDEs, then
each weak convergence is for a fixed $x$, and the choice of
subsequence may depend on $x$. However, this will cause serious
problems when one considers weak solutions. Our approach to avoid
this essential difficulty is to find a subsequence of the weak
convergence in the space $M^2([t,T],
L_{\rho}^2({\mathbb{R}^{d}};{\mathbb{R}^{1}}))\times
M^{2}([t,T];L_{\rho}^2({\mathbb{R}^{d}};{\mathbb{R}^{d}}))$. The
whole point and major difficulty of this approach are to pass the
limit term by term in the approximating equation to the desired
limit. This is achieved in our paper by obtaining a strong
convergent subsequence of $(Y_s^{t,x,n}, Z_{s}^{t,x,n})$ in
$M^2([t,T], L_{\rho}^2({\mathbb{R}^{d}};{\mathbb{R}^{1}}))\times
M^{2}([t,T];L_{\rho}^2({\mathbb{R}^{d}};{\mathbb{R}^{d}}))$ using
the Rellich-Kondrachov compactness theorem and generalized
equivalence of norm principle as we have already mentioned.

\section{The main results}
\setcounter{equation}{0}

In this paper, we study the weak solutions of a class of parabolic
PDEs with p-growth coefficients, their corresponding backward
stochastic differential equations (BSDEs) in a Hilbert space
($\rho$-weighted $L^2$ space) and the probabilistic representation
of the weak solutions of (\ref{1jan2009}) by using the solutions of
BSDEs. For an arbitrary fixed $t\in[0,T]$, we start from the following SDE:
\begin{eqnarray}\label{a}
X_s^{t,x}=x+\int_t^sb(X_r^{t,x})dr+\int_t^s\sigma(X_r^{t,x})dW_r,\ \
\ s\geq t,
\end{eqnarray}
where $W$ is a $\mathbb{R}^d$ Brownian motion on a probability space
$(\Omega, \mathscr{F}, P)$, and
$b:\mathbb{R}^d\rightarrow\mathbb{R}^d$,
$\sigma:\mathbb{R}^d\rightarrow\mathbb{R}^{d\times d}$ are
measurable. For given $t\in[0,T]$ in (\ref{a}), we consider a slightly more general BSDEs for $s\geq t$ by allowing
$f$ depending on time explicitly:
\begin{eqnarray}\label{t}
Y_s^{t,x}=h(X_T^{t,x})+\int_s^Tf(r,X_r^{t,x},Y_r^{t,x},Z_r^{t,x})dr-\int_s^T\langle
Z_r^{t,x},dW_r\rangle,
\end{eqnarray}
where
$f:[0,T]\times\mathbb{R}^d\times\mathbb{R}^1\times\mathbb{R}^d\rightarrow\mathbb{R}^1$
and $h:\mathbb{R}^d\rightarrow\mathbb{R}^1$ are measurable. More
conditions on $b$, $\sigma$, $f$ are needed and will be specified
later. The Hilbert space
$L_{\rho}^2({\mathbb{R}^{d}};{\mathbb{R}^{k}})$ is the space
containing all Borel measurable functions $l$:
$\mathbb{R}^d\to\mathbb{R}^k$ such that $\int_{\mathbb{R}^{d}}<l(x),
l(x)>\rho^{-1}(x)dx<\infty$, with the inner product
\begin{eqnarray*}
\langle
u_1,u_2\rangle=\int_{\mathbb{R}^d}<u_1(x),u_2(x)>\rho^{-1}(x)dx,
\end{eqnarray*}
where $\rho(x)=(1+|x|)^q$, $q>d$, is a weight function. 
The Banach space $L_{\rho}^{2p}({\mathbb{R}^{d}};{\mathbb{R}^{1}})$
is the space containing all Borel measurable functions $l$:
$\mathbb{R}^d\to\mathbb{R}^1$ such that
$\int_{\mathbb{R}^{d}}l^{2p}(x)\rho^{-1}(x)dx<\infty$ with the norm
$||l||_{L_{\rho}^{2p}(R^d)}=(\int_{\mathbb{R}^{d}}l^{2p}(x)\rho^{-1}(x)dx)^{1\over
2p}$. It is easy to see that
$\rho(x):\mathbb{R}^d\longrightarrow\mathbb{R}^1$ is a continuous
positive function satisfying
$\int_{\mathbb{R}^{d}}\rho^{-1}(x)dx<\infty$. 
Note that we can consider more general $\rho$ which satisfies the
above condition and conditions in \cite{ba-ma} and all
the results of this paper still hold. For $k\geq0$, we denote by $C_{b}^k$
\ the set of $C^k$-functions whose partial derivatives of order less
than or equal to $k$ are bounded, 
and denote by
$H^1_\rho$ 
the $\rho$-weighted Sobolev space, i.e. the completion of $C_c^\infty(\mathbb{R}^d;\mathbb{R}^1)$ w.r.t. the norm
$\|\varphi\|^2_{H_{\rho}^{1}(\mathbb{R}^d;\mathbb{R}^1)}=\int_{\mathbb{R}^d}(|\varphi(x)|^2+|\nabla\varphi(x)|^2)\rho^{-1}(x)dx$. Now we
assume the following conditions for the coefficients in SDE
(\ref{a}) and BSDE (\ref{t}):
\begin{description}
\item[(H.1).] For a given $p\geq
1$, $\int_{\mathbb{R}^d}|h(x)|^{2p}\rho^{-1}(x)dx<\infty$.
\item[(H.2).] There exists a constant $C\geq0$ and a function ${f}_0$ with $\int_0^T\int_{\mathbb{R}^d}|{f}_0(s,x)|^{2p}\rho^{-1}(x)dxds<\infty$ s.t. $|f(s,x,y,z)|\leq C(|{f}_0(s,x)|+|y|^p+|z|)$, where $p$ is
the same as in (H.1).
\item[(H.3).] There exists a
constant $\mu\in\mathbb{R}^{1}$ s.t. for any $s\in[0,T]$, $y_1,
y_2\in {\mathbb{R}^{1}}$, $x,z\in{\mathbb{R}^{d}}$,
\begin{eqnarray*}
(y_1-y_2)\big(f(s,x,y_1,z)-f(s,x,y_2,z)\big)\leq\mu{|y_1-y_2|}^2.
\end{eqnarray*}
\item[(H.4).] The function $(y,z)\rightarrow f(s,x,y,z)$ is continuous and $z\rightarrow f(s,x,y,z)$ is globally
Lipschitz continuous with Lipschitz constant $L\geq0$, i.e. for any
$s\in[0,T]$, $y\in {\mathbb{R}^{1}}$,
$x,z_1,z_2\in{\mathbb{R}^{d}}$,
\begin{eqnarray*}
|f(s,x,y,z_1)-f(s,x,y,z_2)|\leq L|z_1-z_2|.
\end{eqnarray*}
\item[(H.5).] The coefficients $b\in C_{b}^2(\mathbb{R}^{d};\mathbb{R}^{d})$, $\sigma\in
C_{b}^3(\mathbb{R}^{d};\mathbb{R}^{d}\times\mathbb{R}^{d})$ and
$\sigma$ satisfies the uniform ellipticity condition, i.e. there
exists a constant $D>0$ s.t. $\xi^*(\sigma\sigma^*)(x)\xi\geq
D\xi^*\xi$ for any $\xi\in\mathbb{R}^{d}$.
\end{description}

Condition (H.5) guarantees the existence of the flow of diffeomorphism. This is essential in the equivalence of norm principle (Lemma \ref{qi045}), which together with the uniform ellipticity condition implies the equivalence of the Sobolev norm of the solution of the PDE and $L^2(\mathbb{R}^d;\mathbb{R}^1)\times L^2(\mathbb{R}^d;\mathbb{R}^d)$ norm of the solution of the BSDE. See (\ref{am}) in Section 4.

It is easy to see that for a.a. $x\in\mathbb{R}^d$,
$(Y_s^{t,x},Z_s^{t,x})$ solves BSDE (\ref{t}) if and only if
$(\tilde{Y}_s^{t,x},\tilde{Z}_s^{t,x})=({\rm e}^{\mu
s}Y_s^{t,x},{\rm e}^{\mu s}Z_s^{t,x})$ solves the following BSDE:
\begin{eqnarray}\label{e}
\tilde{Y}_s^{t,x}={\rm e}^{\mu
T}h(X_T^{t,x})+\int_s^T\tilde{f}(r,X_r^{t,x},\tilde{Y}_r^{t,x},\tilde{Z}_r^{t,x})dr-\int_s^T\langle
\tilde{Z}_r^{t,x},dW_r\rangle,
\end{eqnarray}
where $\tilde{f}(r,x,y,z)={\rm e}^{\mu r}f(r,x,{\rm e}^{-\mu
r}y,{\rm e}^{-\mu r}z)-\mu y$. We can verify that $\tilde{f}$
satisfies Conditions (H.2), (H.3) and (H.4). But, by Condition
(H.3), for $y_1, y_2\in\mathbb{R}^1$, and $x,z\in\mathbb{R}^d$,
\begin{eqnarray*}
&&(y_1-y_2)\big(\tilde{f}(s,x,y_1,z)-\tilde{f}(s,x,y_2,z)\big)\\
&=&{\rm e}^{2\mu s}({\rm e}^{-\mu s}y_1-{\rm e}^{-\mu s}y_2)\big(f(s,x,{\rm e}^{-\mu s}y_1,{\rm e}^{-\mu s}z)-f(s,x,{\rm e}^{-\mu s}y_2,{\rm e}^{-\mu s}z)\big)-\mu(y_1-y_2)(y_1-y_2)\\
&\leq&\mu{\rm e}^{2\mu s}|{\rm e}^{-\mu s}y_1-{\rm e}^{-\mu
s}y_2|^2-\mu|y_1-y_2|^2=0.
\end{eqnarray*}

Now we give the definition for the solution of BSDE (\ref{t}) in the
$\rho$-weighted $L^2$ space. First define the space for the solution
$(Y_\cdot^{t,\cdot},Z_\cdot^{t,\cdot})$. We denote by ${\cal N}$ the
class of $P$-null sets of ${\mathscr{F}}$ and let
$\mathscr{F}_t\triangleq\mathscr{F}_t^W\bigvee{\cal N}$, for $0\leq
t\leq T$. We recall some definitions.
\begin{defi}\ (Definitions 2.2 in \cite{zh-zh1})\label{zhao005}
Let $\mathbb{S}$ be a separable Banach space with norm $\|\cdot\|_\mathbb{S}$
and Borel $\sigma$-field $\mathscr{S}$ and $q\geq2$ be a real
number. We denote by $M^{q}([t,T];\mathbb{S})$ the set of
$\mathscr{B}([t,T])\otimes\mathscr{F}/\mathscr{S}$ measurable random
processes $\{\phi(s)\}_{t\leq s\leq T}$ with values in $\mathbb{S}$
satisfying
\begin{enumerate-roman}
\item $\phi(s):\Omega\rightarrow\mathbb{S}$ is $\mathscr{F}_{s}$ measurable for $t\leq s\leq T$;
\item $E[\int_{t}^{T}\|\phi(s)\|_\mathbb{S}^qds]<\infty$.
\end{enumerate-roman}
Also we denote by $S^{q}([t,T];\mathbb{S})$ the set of
$\mathscr{B}([t,T])\otimes\mathscr{F}/\mathscr{S}$ measurable random
processes $\{\psi(s)\}_{t\leq s\leq T}$ with values in $\mathbb{S}$
satisfying
\begin{enumerate-roman}
\item $\psi(s):\Omega\rightarrow\mathbb{S}$ is
$\mathscr{F}_{s}$ measurable for $t\leq s\leq T$ and
$\psi(\cdot,\omega)$ is continuous $P$-a.s.;
\item $E[\sup_{t\leq s\leq T}\|\psi(s)\|_\mathbb{S}^q]<\infty$.
\end{enumerate-roman}
\end{defi}
\begin{defi}\ (Definitions 3.1 in \cite{zh-zh1})\label{qi051}
A pair of processes $(Y_{s}^{t,x},Z_{s}^{t,x})$ is called a solution
of BSDE (\ref{t}) if $(Y_{\cdot}^{t,\cdot},Z_{\cdot}^{t,\cdot})\in
S^{2p}([t,T];L_{\rho}^{2p}({\mathbb{R}^{d}};{\mathbb{R}^{1}}))\times
M^{2}([t,T];L_{\rho}^2({\mathbb{R}^{d}};{\mathbb{R}^{d}}))$ and
$(Y_s^{t,x},Z_s^{t,x})$ satisfies (\ref{t}) for a.a. $x$, with
probability one.
\end{defi}
\begin{rmk}\label{16}
Due to the density of $C_c^{0}(\mathbb{R}^d;\mathbb{R}^1)$ in
$L_{\rho}^2({\mathbb{R}^{d}};{\mathbb{R}^{d}})$,
$(Y_s^{t,x},Z_s^{t,x})$ satisfies (\ref{t}) for a.a. $x$ with probability one is
equivalent to that for an arbitrary $\varphi\in
C_c^{0}(\mathbb{R}^d;\mathbb{R}^1)$, $({Y}_s^{t,x},{Z}_s^{t,x})$
satisfies
\begin{eqnarray*}
\int_{\mathbb{R}^{d}}Y_s^{t,x}\varphi(x)dx&=&\int_{\mathbb{R}^{d}}h(X_T^{t,x})\varphi(x)dx+\int_{s}^{T}\int_{\mathbb{R}^{d}}f(r,X_{r}^{t,x},Y_{r}^{t,x},Z_{r}^{t,x})\varphi(x)dxdr\nonumber\\
&&-\int_{s}^{T}\langle
\int_{\mathbb{R}^{d}}Z_{r}^{t,x}\varphi(x)dx,dW_r\rangle\ \ \ P-{\rm
a.s.}
\end{eqnarray*}
\end{rmk}

Since $(Y_\cdot^{t,\cdot},Z_\cdot^{t,\cdot})\in
S^{2p}([t,T];L_{\rho}^{2p}({\mathbb{R}^{d}};{\mathbb{R}^{1}}))\times
M^{2}([t,T];L_{\rho}^2({\mathbb{R}^{d}};{\mathbb{R}^{d}}))$ if and
only if $(\tilde{Y}_\cdot^{t,\cdot},\tilde{Z}_\cdot^{t,\cdot})\in
S^{2p}([t,T];L_{\rho}^{2p}({\mathbb{R}^{d}};{\mathbb{R}^{1}}))\times
M^{2}([t,T];L_{\rho}^2({\mathbb{R}^{d}};{\mathbb{R}^{d}}))$, so we
claim $(Y_s^{t,x},Z_s^{t,x})$ is the solution of BSDE (\ref{t}) in
the $\rho$-weighted $L^2$ space if and only if
$(\tilde{Y}_s^{t,x},\tilde{Z}_s^{t,x})$ is the solution of BSDE
(\ref{e}) in $\rho$-weighted $L^2$ space. Therefore we can replace,
without losing any generality, Condition (H.3) by
\begin{description}
\item[${\rm (H.3)}^*$.] For any $s\in[0,T]$, $y_1,
y_2\in {\mathbb{R}^{1}}$, $x,z\in{\mathbb{R}^{d}}$,
\begin{eqnarray*}
(y_1-y_2)\big(f(s,x,y_1,z)-f(s,x,y_2,z)\big)\leq0.
\end{eqnarray*}
\end{description}

The main purpose of this paper is to prove the following two
theorems. The first one is about the existence and uniqueness of
solutions to BSDE (\ref{t}):
\begin{thm}\label{21}
Under Conditions (H.1), (H.2), (H.3)$^*$, (H.4) and (H.5), BSDE
(\ref{t}) has a unique solution $(Y_\cdot^{t,\cdot},
Z_\cdot^{t,\cdot})\in
S^{2p}([t,T];L_{\rho}^{2p}({\mathbb{R}^{d}};{\mathbb{R}^{1}}))\times
M^{2}([t,T];L_{\rho}^2({\mathbb{R}^{d}};{\mathbb{R}^{d}}))$.
\end{thm}

We will establish a connection between BSDE (\ref{t}) and the following PDE with p-growth
coefficients:
\begin{numcases}{}\label{z}
{{\partial u}\over{\partial
t}}(t,x)=-\mathscr{L}u(t,x)-f\big(t,x,u(t,x),(\sigma^*\nabla
u)(s,x)\big),\ \ \ \ 0\leq t\leq
T,\nonumber\\
u(T,x)=h(x).
\end{numcases}


Noticing $f$ is of p-growth on $y$, we recall the definition for the
weak solution of PDE (\ref{z}):
\begin{defi}\label{12}
Function $u$ is called the weak solution of PDE (\ref{z}) if
$(u,\sigma^*\nabla u)\in
L^{2p}([0,T];L_{\rho}^{2p}({\mathbb{R}^{d}};{\mathbb{R}^{1}}))\\\times
L^{2}([0,T];L_{\rho}^2({\mathbb{R}^{d}};{\mathbb{R}^{d}}))$ and for
an arbitrary $\varphi\in C_c^{\infty}(\mathbb{R}^d;\mathbb{R}^1)$,
\begin{eqnarray}\label{ae}
&&\int_{\mathbb{R}^{d}}u(t,x)\varphi(x)dx-\int_{\mathbb{R}^{d}}u(T,x)\varphi(x)dx-{1\over2}\int_{t}^{T}\int_{\mathbb{R}^{d}}\big((\sigma^*\nabla u)(s,x)\big)^*(\sigma^*\nabla\varphi)(x)dxds\nonumber\\
&&-\int_{t}^{T}\int_{\mathbb{R}^{d}}u(s,x)div\big((b-\tilde{A})\varphi\big)(x)dxds\nonumber\\
&=&\int_{t}^{T}\int_{\mathbb{R}^{d}}f\big(s,x,u(s,x),(\sigma^*\nabla
u)(s,x)\big)\varphi(x)dxds,\ \ t\in[0,T].
\end{eqnarray}
Here $\tilde{A}_j\triangleq{1\over2}\sum_{i=1}^d{\partial
a_{ij}(x)\over\partial x_i}$, and
$\tilde{A}=(\tilde{A}_1,\tilde{A}_2,\cdot\cdot\cdot,\tilde{A}_d)^*$.
\end{defi}

The other main theorem is the probabilistic representation of PDE
(\ref{z}) in the $\rho$-weighted $L^2$ space through its
corresponding BSDE:
\begin{thm}\label{20}
Define $u(t,x)=Y_t^{t,x}$, where $(Y_s^{t,x},Z_s^{t,x})$ is the
solution of BSDE (\ref{t}) under Conditions (H.1), (H.2), (H.3)$^*$,
(H.4) and (H.5), then $u(t,x)$ is the unique weak solution of PDE
(\ref{z}). Moreover, let $u$ be a representative in the equivalence class of the solution of the PDE (\ref{z}) in $L^{2p}([0,T];L_{\rho}^{2p}({\mathbb{R}^{d}};{\mathbb{R}^{1}}))$ with $\sigma^*\nabla u\in
L^{2}([0,T];L_{\rho}^2({\mathbb{R}^{d}};{\mathbb{R}^{d}}))$, then $u(t,x)=Y_t^{t,x}$ for a.a. $t\in[0,T]$, a.a. $x\in\mathbb{R}^d$ and
\begin{eqnarray}\label{ce}
u(s,X_s^{t,x})=Y_s^{t,x},\ (\sigma^*\nabla u)(s,X^{t,x}_s)=Z_s^{t,x}
\ {\rm for}\ {\rm a.a.}\ s\in[t,T],\ {\rm a.a.}\ x\in\mathbb{R}^{d}\ {\rm a.s.}
\end{eqnarray}
\end{thm}
We give the proofs of these two theorems in the latter sections.


In Sections 3-5, by making use of truncated BSDEs, we first deal
with BSDE (\ref{t}). To prove BSDE (\ref{t}) has a unique solution,
we use the Alaoglu lemma to derive a weakly convergent sequence in
Section 3 and further use the equivalence of norm principle and
Rellich-Kondrachov Compactness Theorem to get a strongly convergent
sequence in Section 4. Then we complete the proofs of Theorem
\ref{21} in Section 5 and consider the corresponding PDE (\ref{z})
to obtain Theorem \ref{20} in Section 6 which gives the
probabilistic representation to the weak solution of PDE (\ref{z}).

\begin{rmk}
Let $u$ be the weak solution of PDE (\ref{z}) with coefficient
$f\big(x,u,(\sigma^*\nabla u)\big)$ which is independent of $t$, we
can see easily that $v(t)\triangleq u(T-t)$ is the unique weak
solution of PDE (\ref{1jan2009}).
\end{rmk}

\section{The weak convergence}
\setcounter{equation}{0}

Assume $f$ satisfies Conditions (H.2), ${\rm (H.3)}^*$ and (H.4). We
first use a standard cut-off technique to study a sequence of BSDEs
with nonlinear function $f_n$ satisfying the linear growth condition
on $y$. The
$S^{2p}([t,T];L_{\rho}^{2p}({\mathbb{R}^{d}};{\mathbb{R}^{1}}))\times
M^{2}([t,T];L_{\rho}^2({\mathbb{R}^{d}};{\mathbb{R}^{d}}))$ valued
solution for this kind equation was studied in \cite{zh-zh2}. For
this, we define for each $n\in N$
\begin{eqnarray}\label{13}
f_n(s,x,y,z)=f\big(s,x,\Pi_n(y),z\big),
\end{eqnarray}
where $\Pi_n(y)={\inf(n,|y|)\over|y|}y$. Then
$f_n:[0,T]\times\mathbb{R}^d\times\mathbb{R}^1\times\mathbb{R}^d\rightarrow\mathbb{R}^1$
satisfies
\begin{description}
\item[(H.2)$'$.] For any $s\in[0,T]$, $y\in
{\mathbb{R}^{1}}$, $x,z\in{\mathbb{R}^{d}}$ and the constant $C$
given in (H.2), $$|f_n(s,x,y,z)|\leq C(|f_0(s,x)|+|n|^p+|z|).$$
\item[(H.3)$'$.] For any $s\in[0,T]$, $y_1,
y_2\in {\mathbb{R}^{1}}$, $x\in{\mathbb{R}^{d}}$,
\begin{eqnarray*}
(y_1-y_2)\big(f_n(s,x,y_1,z)-f_n(s,x,y_2,z)\big)\leq0.
\end{eqnarray*}
\item[(H.4)$'$.] The function $(y,z)\rightarrow f_n(s,x,y,z)$ is
 continuous, and for any $s\in[0,T]$, $y\in {\mathbb{R}^{1}}$,
$x,z_1,z_2\in{\mathbb{R}^{d}}$ and the constant $L$ given in (H.4),
\begin{eqnarray*}
|f_n(s,x,y,z_1)-f_n(s,x,y,z_2)|\leq L|z_1-z_2|.
\end{eqnarray*}
\end{description}
To see (H.3)$'$, if $\Pi_n(y_1)=\Pi_n(y_2)$, it is obvious; if
$\Pi_n(y_1)\neq\Pi_n(y_2)$, then
\begin{eqnarray*}
&&(y_1-y_2)\big(f_n(s,x,y_1,z)-f_n(s,x,y_2,z)\big)\nonumber\\
&=&(\Pi_n(y_1)-\Pi_n(y_2))\big(f(s,x,\Pi_n(y_1),z)-f(s,x,\Pi_n(y_2),z)\big){{y_1-y_2}\over{\Pi_n(y_1)-\Pi_n(y_2)}}\leq0.
\end{eqnarray*}

We then study the following BSDE with the global Lipschitz
coefficient $f_n$:
\begin{eqnarray}\label{c}
Y_s^{t,x,n}=h(X_T^{t,x})+\int_s^Tf_n(r,X_r^{t,x},Y_r^{t,x,n},Z_r^{t,x,n})dr-\int_s^T\langle
Z_r^{t,x,n},dW_r\rangle.
\end{eqnarray}
Notice that under the conditions of Theorem \ref{21}, the
coefficients $h$ and $f_n$ satisfy Conditions (H.1), (H.2)$'$ and
(H.4)$'$. Hence by Theorems 2.2 and 2.3 in \cite{zh-zh2}, we have
the following proposition:
\begin{prop}\label{7} (\cite{zh-zh2})
Under the conditions of Theorem \ref{21}, for $f_n$ defined in
(\ref{13}), BSDE (\ref{c}) has a unique solution $(Y_s^{t,x,n},
Z_s^{t,x,n})\in
S^{2}([t,T];L_{\rho}^2({\mathbb{R}^{d}};{\mathbb{R}^{1}}))\times
M^{2}([t,T];L_{\rho}^2({\mathbb{R}^{d}};{\mathbb{R}^{d}}))$. If we
define $Y_t^{t,x,n}=u_n(t,x)$, then $u_n(t,x)$ is the unique weak
solution of the following PDE
\begin{numcases}{}\label{f}
{{\partial u_n}\over{\partial
t}}(t,x)=-\mathscr{L}u_n(t,x)-f_n\big(t,x,u_n(t,x),(\sigma^*\nabla
u)(t,x)\big),\ \ \ \ 0\leq t\leq
T,\nonumber\\
u_n(T,x)=h(x).
\end{numcases}
Moreover,
\begin{eqnarray}\label{zz16}
u_n(s,X_s^{t,x})=Y_s^{t,x,n},\ (\sigma^*\nabla
u_n)(s,X^{t,x}_s)=Z_s^{t,x,n}\ {\rm for}\ {\rm a.a.}\ s\in[t,T],\ {\rm a.a.}\ x\in\mathbb{R}^{d}\ {\rm a.s.}
\end{eqnarray}
\end{prop}

The key is to pass the limits in (\ref{c}) and (\ref{f}) in some
desired sense. For this we need some estimates that go beyond those
in \cite{zh-zh1} and \cite{zh-zh2}. Before we derive some useful
estimations to the solution of BSDEs (\ref{c}), we give the
generalized equivalence of norm principle which is an extension of
equivalence of norm principle given in \cite {ku1}, \cite{ba-le},
\cite{ba-ma} to the cases when $\varphi$ and $\Psi$ are random.
\begin{lem}\label{qi045}(generalized equivalence
of norm principle \cite{zh-zh1}) Let $\rho$ be the weight function
defined at the beginning of Section 1 and $X$ be a diffusion process
defined in (\ref{a}), where the coefficients $b\in C_{b}^2(\mathbb{R}^{d};\mathbb{R}^{d})$, $\sigma\in
C_{b}^3(\mathbb{R}^{d};\mathbb{R}^{d}\times\mathbb{R}^{d})$. If $s\in[t,T]$,
$\varphi:\Omega\times\mathbb{R}^d\rightarrow\mathbb{R}^1$ is
independent of the $\sigma$-field $\mathscr{F}^W_{t,s}\triangleq\sigma\{W_{r}-W_{t},\ t\leq r\leq
s\}$ and $\varphi\rho^{-1}\in L^1(\Omega\times\mathbb{R}^{d})$,
then there exist two constants $c>0$ and $C>0$ such that
\begin{eqnarray*}
cE[\int_{\mathbb{R}^{d}}|\varphi(x)|\rho^{-1}(x)dx]\leq
E[\int_{\mathbb{R}^{d}}|\varphi(X_{s}^{t,x})|\rho^{-1}(x)dx]\leq
CE[\int_{\mathbb{R}^{d}}|\varphi(x)|\rho^{-1}(x)dx].
\end{eqnarray*}
Moreover if
$\Psi:\Omega\times[t,T]\times\mathbb{R}^d\rightarrow\mathbb{R}^1$,
$\Psi(s,\cdot)$ is independent of $\mathscr{F}^W_{t,s}$ and
$\Psi\rho^{-1}\in L^1(\Omega\times[t,T]\times\mathbb{R}^{d})$,
then
\begin{eqnarray*}
&&cE[\int_{t}^{T}\int_{\mathbb{R}^{d}}|\Psi(s,x)|\rho^{-1}(x)dxds]\leq E[\int_{t}^{T}\int_{\mathbb{R}^{d}}|\Psi(s,X_{s}^{t,x})|\rho^{-1}(x)dxds]\\
&\leq&CE[\int_{t}^{T}\int_{\mathbb{R}^{d}}|\Psi(s,x)|\rho^{-1}(x)dxds].
\end{eqnarray*}
\end{lem}

First we deduce a useful estimate.
\begin{lem}\label{1}
Under Conditions (H.1), (H.2), (H.3)$^*$, (H.4) and (H.5), if
$(Y_\cdot^{t,\cdot,n},Z_\cdot^{t,\cdot,n})$ is the solution of BSDE
(\ref{c}), then we have
\begin{eqnarray*}
E[\int_t^T\sup_n\int_{\mathbb{R}^d}|Y_s^{t,x,n}|^{2p}\rho^{-1}(x)dxds]+\sup_nE[\int_{t}^{T}\int_{\mathbb{R}^{d}}{{|{Y}_s^{t,x,n}|}^{2p-2}}|{Z}_{s}^{t,x,n}|^2\rho^{-1}(x)dxds]<\infty.
\end{eqnarray*}
\end{lem}
{\em Proof}. For $M$, $N>0$ and $m\geq2$, define
\begin{eqnarray*}
\psi_M(y)=y^2I_{\{-M\leq y<M\}}+M(2y-M)I_{\{y\geq
M\}}-M(2y+M)I_{\{y<-M\}}
\end{eqnarray*}
and
\begin{eqnarray*}
\varphi_{N,m}(y)=y^{m\over2}I_{\{0\leq
y<N\}}+N^{{m-2}\over2}({m\over2}y-{{m-2}\over2}N)I_{\{y\geq N\}}.
\end{eqnarray*}
Applying It$\hat {\rm o}$'s formula to ${\rm
e}^{Kr}\varphi_{N,m}\big(\psi_M(Y_{r}^{t,x,n})\big)$ for a.a.
$x\in\mathbb{R}^d$, we have
\begin{eqnarray}\label{af}
&&{\rm e}^{Ks}\varphi_{N,m}\big(\psi_M(Y_{s}^{t,x,n})\big)+K\int_{s}^{T}{\rm e}^{Kr}\varphi_{N,m}\big(\psi_M(Y_{r}^{t,x,n})\big)dr\nonumber\\
&&+{1\over2}\int_{s}^{T}{\rm e}^{Kr}\varphi^{''}_{N,m}\big(\psi_M(Y_{r}^{t,x,n})\big)|\psi_M^{'}(Y_{r}^{t,x,n})|^2|Z_{r}^{t,x,n}|^2dr\nonumber\\
&&+\int_{s}^{T}{\rm e}^{Kr}\varphi^{'}_{N,m}\big(\psi_M(Y_{r}^{t,x,n})\big)I_{\{-M\leq{Y}_{r}^{t,x,n}<M\}}|{Z}_{r}^{t,x,n}|^2dr\nonumber\\
&=&{\rm e}^{KT}\varphi_{N,m}\big(\psi_M(h(X_{T}^{t,x}))\big)+\int_{s}^{T}{\rm e}^{Kr}\varphi^{'}_{N,m}\big(\psi_M(Y_{r}^{t,x,n})\big)\psi_M^{'}(Y_{r}^{t,x,n}){f_n}(r,X_{r}^{t,x},Y_{r}^{t,x,n},Z_{r}^{t,x,n})dr\nonumber\\
&&-\int_{s}^{T}\langle{\rm
e}^{Kr}\varphi^{'}_{N,m}\big(\psi_M(Y_{r}^{t,x,n})\big)\psi_M^{'}(Y_{r}^{t,x,n}){Z}_{r}^{t,x,n},dW_r\rangle.
\end{eqnarray}
From \cite{zh-zh1}, we note first
$(Y_{\cdot}^{t,\cdot,n},Z_{\cdot}^{t,\cdot,n})\in
S^{2}([0,T];L_{\rho}^2({\mathbb{R}^{d}};{\mathbb{R}^{1}}))\times
M^{2}([0,T];L_{\rho}^2({\mathbb{R}^{d}};{\mathbb{R}^{d}}))$. Also it
is obvious that
$\varphi^{'}_{N,m}\big(\psi_M(Y_{r}^{t,x,n})\big)\psi_M^{'}(Y_{r}^{t,x,n})$
is bounded, hence we can use the stochastic Fubini theorem and take
the conditional expectation w.r.t. $\mathscr{F}_s$. Note that the
stochastic integral has zero conditional expectation. So if we
define ${\psi^{'}_M(y)\over y}=2$ when $y=0$, we have
\begin{eqnarray*}
&&\int_{\mathbb{R}^d}{\rm e}^{Ks}\varphi_{N,m}\big(\psi_M(Y_{s}^{t,x,n})\big)\rho^{-1}(x)dx+E[K\int_{s}^{T}\int_{\mathbb{R}^d}{\rm e}^{Kr}\varphi_{N,m}\big(\psi_M(Y_{r}^{t,x,n})\big)\rho^{-1}(x)dxdr|\mathscr{F}_s]\nonumber\\
&&+{1\over2}E[\int_{s}^{T}\int_{\mathbb{R}^d}{\rm e}^{Kr}\varphi^{''}_{N,m}\big(\psi_M(Y_{r}^{t,x,n})\big)|\psi_M^{'}(Y_{r}^{t,x,n})|^2|Z_{r}^{t,x,n}|^2\rho^{-1}(x)dxdr|\mathscr{F}_s]\nonumber\\
&&+E[\int_{s}^{T}\int_{\mathbb{R}^d}{\rm e}^{Kr}\varphi^{'}_{N,m}\big(\psi_M(Y_{r}^{t,x,n})\big)I_{\{-M\leq{Y}_{r}^{t,x,n}<M\}}|{Z}_{r}^{t,x,n}|^2\rho^{-1}(x)dxdr|\mathscr{F}_s]\nonumber\\
&=&E[\int_{\mathbb{R}^d}{\rm e}^{KT}\varphi_{N,m}\big(\psi_M(h(X_{T}^{t,x}))\big)\rho^{-1}(x)dx|\mathscr{F}_s]\nonumber\\
&&+E[\int_{s}^{T}\int_{\mathbb{R}^d}{\rm e}^{Kr}\varphi^{'}_{N,m}\big(\psi_M(Y_{r}^{t,x,n})\big)\psi_M^{'}(Y_{r}^{t,x,n}){f_n}(r,X_{r}^{t,x},Y_{r}^{t,x,n},Z_{r}^{t,x,n})\rho^{-1}(x)dxdr|\mathscr{F}_s]\nonumber\\
&=&E[\int_{\mathbb{R}^d}{\rm e}^{KT}\varphi_{N,m}\big(\psi_M(h(X_{T}^{t,x}))\big)\rho^{-1}(x)dx|\mathscr{F}_s]\nonumber\\
&&+E[\int_{s}^{T}\int_{\mathbb{R}^d}{\rm e}^{Kr}\varphi^{'}_{N,m}\big(\psi_M(Y_{r}^{t,x,n})\big){\psi_M^{'}(Y_{r}^{t,x,n})\over{Y_{r}^{t,x,n}}}Y_{r}^{t,x,n}\nonumber\\
&&\ \ \ \ \ \ \ \ \ \ \ \ \ \ \times\big({f}_n(r,X_{r}^{t,x},Y_{r}^{t,x,n},Z_{r}^{t,x,n})-f_n(r,X_{r}^{t,x},0,Z_{r}^{t,x,n})\big)\rho^{-1}(x)dxdr|\mathscr{F}_s]\nonumber\\
&&+E[\int_{s}^{T}\int_{\mathbb{R}^d}{\rm e}^{Kr}\varphi^{'}_{N,m}\big(\psi_M(Y_{r}^{t,x,n})\big)\psi_M^{'}(Y_{r}^{t,x,n})\nonumber\\
&&\ \ \ \ \ \ \ \ \ \ \ \ \ \times\big({f}_n(r,X_{r}^{t,x},0,Z_{r}^{t,x,n})-f_n(r,X_{r}^{t,x},0,0)\big)\rho^{-1}(x)dxdr|\mathscr{F}_s]\nonumber\\
&&+E[\int_{s}^{T}\int_{\mathbb{R}^d}{\rm e}^{Kr}\varphi^{'}_{N,m}\big(\psi_M(Y_{r}^{t,x,n})\big)\psi_M^{'}(Y_{r}^{t,x,n})f_n(r,X_{r}^{t,x},0,0)\rho^{-1}(x)dxdr|\mathscr{F}_s]\nonumber\\
&\leq&E[\int_{\mathbb{R}^d}{\rm e}^{KT}\varphi_{N,m}\big(\psi_M(h(X_{T}^{t,x}))\big)\rho^{-1}(x)dx|\mathscr{F}_s]\nonumber\\
&&+LE[\int_{s}^{T}\int_{\mathbb{R}^d}{\rm
e}^{Kr}|\varphi^{'}_{N,m}\big(\psi_M(Y_{r}^{t,x,n})\big)||\psi_M^{'}(Y_{r}^{t,x,n})||Z_{r}^{t,x,n}|\rho^{-1}(x)dxdr|\mathscr{F}_s]\nonumber\\
&&+E[\int_{s}^{T}\int_{\mathbb{R}^d}{\rm
e}^{Kr}|\varphi^{'}_{N,m}\big(\psi_M(Y_{r}^{t,x,n})\big)||\psi_M^{'}(Y_{r}^{t,x,n})||f(r,X_{r}^{t,x},0,0)|\rho^{-1}(x)dxdr|\mathscr{F}_s].
\end{eqnarray*}
Taking the limit as $M\to \infty$ first, then the limit as $N\to
\infty$, by the monotone convergence theorem and Young inequality,
we have
\begin{eqnarray}\label{cd}
&&\int_{\mathbb{R}^d}{\rm e}^{Ks}|Y_{s}^{t,x,n}|^m\rho^{-1}(x)dx+KE[\int_{s}^{T}\int_{\mathbb{R}^{d}}{\rm e}^{Kr}{|{Y}_r^{t,x,n}|}^m\rho^{-1}(x)dxdr|\mathscr{F}_s]\nonumber\\
&&+{m(m-1)\over2}E[\int_{s}^{T}\int_{\mathbb{R}^{d}}{\rm e}^{Kr}{{|{Y}_r^{t,x,n}|}^{m-2}}|{Z}_{r}^{t,x,n}|^2\rho^{-1}(x)dxdr|\mathscr{F}_s]\nonumber\\
&\leq&E[\int_{\mathbb{R}^d}{\rm e}^{KT}|h(X_{T}^{t,x})|^m\rho^{-1}(x)dx|\mathscr{F}_s]\nonumber\\
&&+mLE[\int_{s}^{T}\int_{\mathbb{R}^d}{\rm e}^{Kr}|Y_{r}^{t,x,n}|^{m-2}|Y_{r}^{t,x,n}||Z_{r}^{t,x,n}|\rho^{-1}(x)dxdr|\mathscr{F}_s]\nonumber\\
&&+mE[\int_{s}^{T}\int_{\mathbb{R}^d}{\rm e}^{Kr}|Y_{r}^{t,x,n}|^{m-2}|Y_{r}^{t,x,n}||f(r,X_{r}^{t,x},0,0)|\rho^{-1}(x)dxdr|\mathscr{F}_s]\nonumber\\
&\leq&E[\int_{\mathbb{R}^d}{\rm e}^{KT}|h(X_{T}^{t,x})|^m\rho^{-1}(x)dx|\mathscr{F}_s]+mL^2E[\int_{s}^{T}\int_{\mathbb{R}^d}{\rm e}^{Kr}|Y_{r}^{t,x,n}|^{m}\rho^{-1}(x)dxdr|\mathscr{F}_s]\nonumber\\
&&+{m\over4}E[\int_{s}^{T}\int_{\mathbb{R}^d}{\rm e}^{Kr}|Y_{r}^{t,x,n}|^{m-2}|Z_{r}^{t,x,n}|^2\rho^{-1}(x)dxdr|\mathscr{F}_s]\nonumber\\
&&+mE[\int_{s}^{T}\int_{\mathbb{R}^d}{\rm e}^{Kr}|Y_{r}^{t,x,n}|^{m}\rho^{-1}(x)dxdr|\mathscr{F}_s]\nonumber\\
&&+{m\over4}E[\int_{s}^{T}\int_{\mathbb{R}^d}{\rm e}^{Kr}|Y_{r}^{t,x,n}|^{m-2}|f(r,X_{r}^{t,x},0,0)|^2\rho^{-1}(x)dxdr|\mathscr{F}_s]\nonumber\\
&\leq&E[\int_{\mathbb{R}^d}{\rm e}^{KT}|h(X_{T}^{t,x})|^m\rho^{-1}(x)dx|\mathscr{F}_s]
+m(L^2+1)E[\int_{s}^{T}\int_{\mathbb{R}^d}{\rm e}^{Kr}|Y_{r}^{t,x,n}|^{m}\rho^{-1}(x)dxdr|\mathscr{F}_s]\nonumber\\
&&+{m\over4}E[\int_{s}^{T}\int_{\mathbb{R}^d}{\rm
e}^{Kr}|Y_{r}^{t,x,n}|^{m-2}|Z_{r}^{t,x,n}|^2\rho^{-1}(x)dxdr|\mathscr{F}_s]\nonumber\\
&&+{m\over4}\cdot{{m-2}\over m}E[\int_{s}^{T}\int_{\mathbb{R}^d}({\rm
e}^{{{m-2}\over m}Kr}|Y_{r}^{t,x,n}|^{m-2})^{{m\over{m-2}}}\rho^{-1}(x)dxdr|\mathscr{F}_s]\nonumber\\
&&+{m\over4}\cdot{2\over m}E[\int_{s}^{T}\int_{\mathbb{R}^d}({\rm
e}^{{2\over m}Kr}|f(r,X_{r}^{t,x},0,0)|^2)^{m\over2}\rho^{-1}(x)dxdr|\mathscr{F}_s].
\end{eqnarray}
Here and in the following, $C_p$ is a generic constant. 
Therefore,
taking $K>
m(L^2+1)+{{m-2}\over4}$, we have
\begin{eqnarray*}
&&E[\int_t^T\sup_n\int_{\mathbb{R}^d}|Y_{s}^{t,x,n}|^m\rho^{-1}(x)dxds]+\sup_nE[\int_{t}^{T}\int_{\mathbb{R}^{d}}{{|{Y}_s^{t,x,n}|}^{m-2}}|{Z}_{s}^{t,x,n}|^2\rho^{-1}(x)dxds]\nonumber\\
&\leq&C_pE[\int_{\mathbb{R}^d}|h(X_T^{t,x})|^m\rho^{-1}(x)dx]+C_pE[\int_{t}^{T}\int_{\mathbb{R}^d}|f_0(s,X_s^{t,x})|^m\rho^{-1}(x)dxds]\nonumber\\
&\leq&C_p\int_{\mathbb{R}^d}|h(x)|^m\rho^{-1}(x)dx+C_p\int_{t}^{T}\int_{\mathbb{R}^d}|f_0(s,x)|^m\rho^{-1}(x)dxds<\infty.
\end{eqnarray*}
In particular, taking $m=2p$, then the lemma follows.
$\hfill\diamond$\\

Taking $m=2$ in the proof of Lemma \ref{1}, we know
\begin{eqnarray}\label{an}
E[\int_t^T\sup_n\int_{\mathbb{R}^d}|Y_s^{t,x,n}|^2\rho^{-1}(x)dxds+\sup_nE[\int_t^T\int_{\mathbb{R}^d}|Z_s^{t,x,n}|^2\rho^{-1}(x)dxds]<\infty.
\end{eqnarray}
Also we have
\begin{eqnarray*}
&&\sup_nE[\int_t^T\int_{\mathbb{R}^d}|f_n(s,X_s^{t,x},Y_s^{t,x,n},Z_s^{t,x,n})|^2\rho^{-1}(x)dxds]\nonumber\\
&\leq&\sup_nE[\int_t^T\int_{\mathbb{R}^d}C(|f_0(s,X_s^{t,x})|^2+|Y_s^{t,x,n}|^{2p}+|Z_s^{t,x,n}|^2)\rho^{-1}(x)dxds]<\infty.
\end{eqnarray*}
The last inequality follows from the equivalence of norms principle
and Lemma \ref{1}. Define
$U_s^{t,x,n}=f_n(s,X_s^{t,x},Y_s^{t,x,n},Z_s^{t,x,n})$, $s\geq t$,
then
\begin{eqnarray}\label{p}
\sup_nE[\int_t^T\int_{\mathbb{R}^d}(|Y_s^{t,x,n}|^2+|Z_s^{t,x,n}|^2+|U_s^{t,x,n}|^2)\rho^{-1}(x)dxds]<\infty.
\end{eqnarray}
Therefore by using the Alaoglu lemma, we know that there exists a
subsequence, still denoted by $(Y_s^{t,x,n}, Z_s^{t,x,n},
U_s^{t,x,n})$, s.t. $(Y_s^{t,x,n}, Z_s^{t,x,n}, U_s^{t,x,n})$
converges weakly to the limit $(Y_s^{t,x}, Z_s^{t,x}, U_s^{t,x})$ in
$L^2_\rho(\Omega\times[t,T]\times\mathbb{R}^d;\mathbb{R}^1\times\mathbb{R}^d\times\mathbb{R}^1)$
(or equivalently
$L^2(\Omega\times[t,T];L_\rho^2(\mathbb{R}^d;\mathbb{R}^1)\times
L_\rho^2(\mathbb{R}^d;\mathbb{R}^d)\times
L_\rho^2(\mathbb{R}^d;\mathbb{R}^1))$. Now we take the weak limit in
$L^2_\rho(\Omega\times[t,T]\times\mathbb{R}^d;\mathbb{R}^1)$) to
BSDEs (\ref{c}), we can verify that $(Y_s^{t,x}, Z_s^{t,x},
U_s^{t,x})$ satisfies the following BSDE:
\begin{eqnarray}\label{d}
Y_s^{t,x}=h(X_T^{t,x})+\int_s^TU_r^{t,x}dr-\int_s^T\langle
Z_r^{t,x},dW_r\rangle.
\end{eqnarray}
For this, we will check the weak convergence term by term. The weak
convergence to the first term is deduced by the definition of
$Y_s^{t,x}$. The weak convergence to the second term is trivial
since $h(X_T^{t,x})$ is independent of $n$. We then check the weak
convergence to the last two terms. Let $\eta\in
L^2_\rho(\Omega\times[t,T]\times\mathbb{R}^d;\mathbb{R}^1)$. Then
noticing
$\int_t^T\sup_nE[\int_s^T\int_{\mathbb{R}^d}|U_r^{t,x,n}|^2\rho^{-1}(x)dxdr]ds<\infty$
due to (\ref{p}), by Lebesgue's dominated convergence theorem, we
have
\begin{eqnarray*}
&&|E[\int_t^T\int_{\mathbb{R}^d}\int_s^T(U_r^{t,x,n}-U_r^{t,x})dr\eta(s,x)\rho^{-1}(x)dxds]|\nonumber\\
&=&|E[\int_t^T\int_s^T\int_{\mathbb{R}^d}(U_r^{t,x,n}-U_r^{t,x})\eta(s,x)\rho^{-1}(x)dxdrds]|\nonumber\\
&\leq&\int_t^T|E[\int_s^T\int_{\mathbb{R}^d}(U_r^{t,x,n}-U_r^{t,x})\eta(s,x)\rho^{-1}(x)dxdr]|ds\longrightarrow0,\
\ \ {\rm as}\ n\to\infty.
\end{eqnarray*}
On the other hand we know for fixed $s$ and $x$, $\eta(s,x)\in
L^2(\Omega)$. So there exists $\varphi(s,x,r)$ s.t.
$\eta(s,x)=E[\eta(s,x)]+\int_t^T\langle\varphi(s,x,r),dW_r\rangle$.
It is easy to see that for a.a. $s\in[t,T]$,
$\varphi(s,\cdot,\cdot)\in
L^2(\Omega\times[t,T]\times\mathbb{R}^d;\mathbb{R}^1)$. Noticing
that
$\int_t^T\sup_nE[\int_s^T\int_{\mathbb{R}^d}|Z_r^{t,x,n}|^2\rho^{-1}(x)dxdr]ds<\infty$
due to (\ref{p}) and using Lebesgue's dominated convergence theorem
again, we obtain
\begin{eqnarray*}
&&|E[\int_t^T\int_{\mathbb{R}^d}\int_s^T\langle Z_r^{t,x,n}-Z_r^{t,x},dW_r\rangle\eta(s,x)\rho^{-1}(x)dxds]|\nonumber\\
&=&|\int_t^T\int_{\mathbb{R}^d}E[\int_s^T\langle Z_r^{t,x,n}-Z_r^{t,x},dW_r\rangle(E[\eta(s,x)]+\int_t^T\langle\varphi(s,x,r),dW_r\rangle)]\rho^{-1}(x)dxds|\nonumber\\
&=&|\int_t^T\int_{\mathbb{R}^d}E[\int_s^T\langle Z_r^{t,x,n}-Z_r^{t,x},\varphi(s,x,r)\rangle dr]\rho^{-1}(x)dxds|\nonumber\\
&\leq&\int_t^T|E[\int_s^T\int_{\mathbb{R}^d}\langle
Z_r^{t,x,n}-Z_r^{t,x},\varphi(s,x,r)\rangle
\rho^{-1}(x)dxdr]|ds\longrightarrow0,\ \ \ {\rm as}\ n\to\infty.
\end{eqnarray*}

Needless to say, if we can show BSDE (\ref{c}) is indeed BSDE
(\ref{t}), then we can say $(Y_s^{t,x},Z_s^{t,x})$ is a solution of
BSDE (\ref{t}). The key is to prove that
$U^{t,x}_s=f(s,X_s^{t,x},Y_s^{t,x},Z_s^{t,x})$ for a.a. $s\in[t,T]$,
$x\in\mathbb{R}^d$ a.s. However, the weak convergence of $Y^n$,
$U^n$ and $Z^n$ are not enough to this. The crucial point in this
analysis is to establish the strong convergence of $Y^n$ and $Z^n$,
which will be done in next section.

\section{The strong convergence and the identification of the limiting BSDEs}
\setcounter{equation}{0}

In this section, we will show that the combination of methods of weak convergence
and strong convergence of a subsequence $(Y_s^{t,x,n},Z_s^{t,x,n})$ 
gives an effective way to prove that the limit
$(Y_s^{t,x},Z_s^{t,x})$ satisfies BSDE (\ref{t}). In contrast, the
direct proof that BSDE (\ref{c}) converges strongly to BSDE
(\ref{t}) by using the strongly convergent subsequence
$(Y_s^{t,x,n},Z_s^{t,x,n})$ without the weak convergence argument
will encounter some complications. This is due to that the dominated
convergence theorem does not seem to apply immediately to
BSDE (\ref{c}). We start from an easy lemma.
\begin{lem}\label{10}
Under the conditions of Theorem \ref{21}, if $u_n(t,x)$ is the weak
solution of PDE (\ref{f}), then
$\sup_n\int_0^T\int_{\mathbb{R}^d}|u_n(s,x)|^{2p}\rho^{-1}(x)dxds<\infty$.
Furthermore, 
\begin{eqnarray*}
\lim_{N\to\infty}\sup_n\int_0^T\int_{\mathbb{R}^d}|u_n(s,x)|^2I_{{U_N}^c}(x)\rho^{-1}(x)dxds=0,
\end{eqnarray*}
where ${U_N}^c=\{x\in\mathbb{R}^d:\ |x|>N\}$. 
\end{lem}
{\em Proof}. By the equivalence of norm principle, (\ref{zz16}) and Lemma \ref{1}, we deduce the $L_\rho^{2p}$ integrability of $u_n$ as follows:
\begin{eqnarray*}
\sup_n\int_0^T\int_{\mathbb{R}^d}|u_n(s,x)|^{2p}\rho^{-1}(x)dxds&\leq&C_p\sup_nE[\int_0^T\int_{\mathbb{R}^d}|u_n(s,X_s^{0,x})|^{2p}\rho^{-1}(x)dxds]\nonumber\\
&=&C_p\sup_nE[\int_0^T\int_{\mathbb{R}^d}|Y_s^{0,x,n}|^{2p}\rho^{-1}(x)dxds]<\infty.
\end{eqnarray*}
Let's then prove the second part of this lemma. Since
$\int_{\mathbb{R}^d}\rho^{-1}(x)dx<\infty$,
\begin{eqnarray*}
&&\lim_{N\to\infty}\sup_n\int_0^T\int_{\mathbb{R}^d}|u_n(s,x)|^2I_{{U_N}^c}(x)\rho^{-1}(x)dxds\nonumber\\
&\leq&\lim_{N\to\infty}\big(\sup_n\int_0^T\int_{\mathbb{R}^d}|u_n(s,x)|^{2p}\rho^{-1}(x)dxds\big)^{1\over p}\big(\int_0^T\int_{\mathbb{R}^d}|I_{{U_N}^c}(x)|^{p\over{p-1}}\rho^{-1}(x)dxds\big)^{{p-1}\over p}\nonumber\\
&\leq&\lim_{N\to\infty}C_p\big(\int_{\mathbb{R}^d}I_{{U_N}^c}(x)\rho^{-1}(x)dx\big)^{{p-1}\over
p}=0.\nonumber
\end{eqnarray*}
$\hfill\diamond$\\

The following two theorems quoted in \cite{ro} will be used in this
section.
\begin{thm}\label{6} (c.f. \cite{ro}) Let $X\subset\subset H\subset Y$ be Banach spaces, with $X$ reflexive. Here $X\subset\subset H$ means
$X$ is compactly embedded in $H$. Suppose that $u_n$ is a sequence
that is uniformly bounded in $L^2([0,T];X)$, and ${du_n/dt}$ is
uniformly bounded in $L^p(0,T;Y)$, for some $p>1$. Then there is a
subsequence that converges strongly in $L^2([0,T];H)$.
\end{thm}
\begin{thm}\label{8} (Rellich-Kondrachov Compactness Theorem c.f. \cite{ro}) Let $B$ be a bounded
 $C^1$ domain in $\mathbb{R}^d$. Then $H^1(B)$ is compactly embedded in $L^2(B)$.
\end{thm}


\begin{lem}\label{5}
Under the conditions of Theorem \ref{21}, if $(Y_s^{t,x,n},
Z_s^{t,x,n})$ is the solution of BSDEs (\ref{c}) and $Y_s^{t,x}$ is
the weak limit of $Y_s^{t,x,n}$ in
$L^2_\rho(\Omega\times[t,T]\times\mathbb{R}^d;\mathbb{R}^1)$, then
there is a subsequence of $Y_s^{t,x,n}$, still denoted by
$Y_s^{t,x,n}$, converging strongly to $Y_s^{t,x}$ in
$L^2(\Omega\times[t,T];L_\rho^2(\mathbb{R}^d;\mathbb{R}^1))$.
\end{lem}
{\em Proof}. Let $u_n(s,x)=Y_s^{s,x,n}$. Then by Proposition
\ref{7}, $u_n(s,X_s^{t,x})=Y_s^{t,x,n}$, $(\sigma^*\nabla
u_n)(s,X^{t,x}_s)=Z_s^{t,x,n}$ for a.a. $s\in[t,T]$,
$x\in\mathbb{R}^{d}$ a.s. Moreover, $u_n(s,x)$ is a weak solution of the PDE (\ref{f}). By the definition of weak solution and the fact that $C_c^\infty(\mathbb{R}^d;\mathbb{R}^1)$ is dense in $H_\rho^1(\mathbb{R}^d;\mathbb{R}^1)$, $u_n(s,x)$ satisfies the
following PDE in ${H_{\rho}^{1}}^{*}(\mathbb{R}^d;\mathbb{R}^1)$:
\begin{eqnarray}\label{g}
du_n(s,x)/ds=-\mathscr{L}u_n(s,x)-f_n\big(s,x,u_n(s,x),(\sigma^*\nabla
u_n)(s,x)\big),\ \ \ \ 0\leq s\leq T.
\end{eqnarray}
To get a strongly convergent subsequence of $Y_s^{t,x,n}$, first note that $u_n$ are uniformly bounded in
$L^2([0,T];H_{\rho}^{1}(\mathbb{R}^d;\mathbb{R}^1))$ by the uniform ellipticity condition
of $\sigma$ and the
equivalence of norm principle:
\begin{eqnarray}\label{am}
&&\sup_n\int_0^T\int_{\mathbb{R}^d}(|u_n(s,x)|^2+|\nabla u_n(s,x)|^2)\rho^{-1}(x)dxds\nonumber\\
&\leq&C_p\sup_n\int_0^T\int_{\mathbb{R}^d}(|u_n(s,x)|^2+|(\sigma^*\nabla u_n)(s,x)|^2)\rho^{-1}(x)dxds\nonumber\\
&\leq&C_p\sup_nE[\int_0^T\int_{\mathbb{R}^d}(|Y_s^{0,x,n}|^2+|Z_s^{0,x,n}|^2)\rho^{-1}(x)dxds]<\infty.
\end{eqnarray}
Then we can deduce that $du_n/ds$ are uniformly bounded in
$L^2([0,T];{H_{\rho}^{1}}^*(\mathbb{R}^d;\mathbb{R}^1))$. For this,
we need to prove that $\mathscr{L}u_n$ and $f_n\in
L^2([0,T];{H_{\rho}^{1}}^*(\mathbb{R}^d;\mathbb{R}^1))$ are
uniformly bounded respectively. First note that for
$i=1,2,\cdots,d$,
\begin{eqnarray*}
|{{\partial \rho^{-1}(x)}\over{\partial
x_i}}|
=|{{-qx_i}\over{(1+|x|)^{q+1}|x|}}|\leq{q\over{(1+|x|)^{q+1}}}\leq
q\rho^{-1}(x).
\end{eqnarray*}
Moreover, recalling the form of $\mathscr{L}$ and noticing the
conditions on $b$ and $\sigma$ in (H.5), we can see that $a_{ij}$
and $b_i$ are uniformly bounded for all $i$, $j$. So for arbitrary
$s\in[0,T]$, $\psi\in C_c^\infty(\mathbb{R}^d;\mathbb{R}^1)$, we
have
\begin{eqnarray*}\label{h}
&&\int_{\mathbb{R}^d}\mathscr{L}u_n(s,x)\cdot\psi(x)\rho^{-1}(x)dx\nonumber\\
&=&\int_{\mathbb{R}^d}\big(-{1\over2}\sum_{i,j=1}^d{{\partial
u_n(s,x)}\over{\partial
x_i}}{\partial(a_{ij}\psi\rho^{-1})(x)\over\partial x_j}-\sum_{i=1}^du_n(s,x){\partial(b_i\psi\rho^{-1})(x)\over{\partial x_i}}\big)dx\nonumber\\
&\leq&\int_{\mathbb{R}^d}(\sum_{i=1}^d|{{\partial
u_n(s,x)}\over{\partial
x_i}}|+|u_n(s,x)|)(\sum_{i,j=1}^d|{\partial(a_{ij}\psi\rho^{-1})(x)\over\partial x_j}|+\sum_{i=1}^d|{\partial(b_i\psi\rho^{-1})(x)\over{\partial x_i}}|)dx\nonumber\\
&\leq&C_p\int_{\mathbb{R}^d}(\sum_{i=1}^d|{{\partial
u_n(s,x)}\over{\partial
x_i}}|+|u_n(s,x)|)(\sum_{j=1}^d|{\partial \psi(x)\over{\partial x_j}}|+|\psi(x)|)\rho^{-1}(x)dx\nonumber\\
&\leq&C_p\sqrt{\int_{\mathbb{R}^d}(\sum_{i=1}^d|{{\partial
u_n(s,x)}\over{\partial
x_i}}|+|u_n(s,x)|)^2\rho^{-1}(x)dx}\sqrt{\int_{\mathbb{R}^d}(\sum_{j=1}^d|{\partial \psi(x)\over{\partial x_j}}|+|\psi(x)|)^2\rho^{-1}(x)dx}\nonumber\\
&\leq&
C_p\|u_n(s,x)\|_{H_{\rho}^{1}(\mathbb{R}^d;\mathbb{R}^1)}\|\psi\|_{{H_{\rho}^{1}}(\mathbb{R}^d;\mathbb{R}^1)}.
\end{eqnarray*}
As $C_c^\infty(\mathbb{R}^d;\mathbb{R}^1)$ is dense in
$H_{\rho}^{1}(\mathbb{R}^d;\mathbb{R}^1)$, therefore for arbitrary
$s\in[0,T]$, it follows that
$\|\mathscr{L}u_n(s,\cdot)\|_{{H_{\rho}^{1}}^*(\mathbb{R}^d;\mathbb{R}^1)}\\\leq
C_p\|u_n(s,\cdot)\|_{H_{\rho}^{1}(\mathbb{R}^d;\mathbb{R}^1)}$ and
by (\ref{am}), we have
\begin{eqnarray*}
\sup_n\|\mathscr{L}u_n\|^2_{L^2([0,T];{H_{\rho}^{1}}^*(\mathbb{R}^d;\mathbb{R}^1))}\leq
C_p\sup_n\int_0^T\int_{\mathbb{R}^d}(|u_n(s,x)|^2+|\nabla
u_n(s,x)|^2)\rho^{-1}(x)dxds<\infty.
\end{eqnarray*}
Also using Lemma \ref{1} and the equivalence of norm principle
again, we obtain
\begin{eqnarray*}
&&\int_0^T\|f_n(s,\cdot,u_n(s,\cdot),(\sigma^*\nabla
u_n)(s,\cdot))\|^2_{L^2_\rho(\mathbb{R}^d;\mathbb{R}^1)}ds\\
&\leq&C_pE[\int_0^T\int_{\mathbb{R}^d}(|f_0(s,x)|^2+|Y_s^{0,x,n}|^{2p}+|Z_s^{0,x,n}|^2)\rho^{-1}(x)dxds]<\infty.
\end{eqnarray*}
Hence $f_n\in
L^2([0,T];{L_{\rho}^{2}}^*(\mathbb{R}^d;\mathbb{R}^1))\subset
L^2([0,T];{H_{\rho}^{1}}^*(\mathbb{R}^d;\mathbb{R}^1))$ and
\begin{eqnarray*}
&&\sup_n\|f_n\|^2_{L^2([0,T];{H_{\rho}^{1}}^*(\mathbb{R}^d;\mathbb{R}^1))}\\
&\leq&C_p\sup_n\int_0^T\|f_n(s,\cdot,u_n(s,\cdot),(\sigma^*\nabla
u_n)(s,\cdot))\|^2_{L^2_\rho(\mathbb{R}^d;\mathbb{R}^1)}ds\nonumber\\
&\leq&C_p\sup_nE[\int_0^T\int_{\mathbb{R}^d}(|f_0(s,x)|^2+|Y_s^{0,x,n}|^{2p}+|Z_s^{0,x,n}|^2)\rho^{-1}(x)dxds]<\infty.
\end{eqnarray*}
Therefore we conclude that $du_n/ds$ are uniformly bounded in
$L^2([0,T];{H_{\rho}^{1}}^*(\mathbb{R}^d;\mathbb{R}^1))$.

Noticing Theorem \ref{8} and applying Theorem \ref{6} with
$X=H_{\rho}^{1}({U_1};\mathbb{R}^1)$,
$H=L_\rho^2({U_1};\mathbb{R}^1)$ and
$Y={H_{\rho}^{1}}^*({U_1};\mathbb{R}^1)$, we are able to extract a
subsequence of $u_n(s,x)$, denoted by $u_{1n}(s,x)$, which converges
strongly in $L^2([0,T];L_\rho^2({U_1};\mathbb{R}^1))$. It is obvious
that this $u_{1n}(s,x)$ satisfies the conditions in Theorem \ref{6}.
Applying Theorem \ref{6} again, we are able to extract a subsequence
of $u_{1n}(s,x)$, denoted by $u_{2n}(s,x)$, that converges strongly
in $L^2([0,T];L_\rho^2({U_2};\mathbb{R}^1))$. Actually we can do
this procedure for all $U_i$, $i=1,2,\cdot\cdot\cdot$. Now we pick
up the diagonal sequence $u_{ii}(s,x)$, $i=1,2,\cdot\cdot\cdot$ and
still denote this sequence by ${u}_n$ for convenience. It is easy to
see that ${u}_n$ converges strongly in all
$L^2([0,T];L_\rho^2({U_i};\mathbb{R}^1))$, $i=1,2,\cdot\cdot\cdot$.
For arbitrary $\varepsilon>0$, noticing Lemma \ref{10}, we can find
$j(\varepsilon)$ large enough such that
\begin{eqnarray*}
\sup_n\int_0^T\int_{{U_{j(\varepsilon)}}^c}
2|{u}_n(s,x)|^2\rho^{-1}(x)dxds<{\varepsilon\over3}.
\end{eqnarray*}
For this $j(\varepsilon)$, there exists $n^*(\varepsilon)>0$ s.t.
when $m,n\geq n^*(\varepsilon)$, we know
\begin{eqnarray*}
\|{u}_m-{u}_n\|^2_{L^2([0,T];L_\rho^2({U_{j(\varepsilon)}};\mathbb{R}^1))}=\int_0^T\int_{U_{j(\varepsilon)}}|{u}_m(s,x)-{u}_n(s,x)|^2\rho^{-1}(x)dxds<{\varepsilon\over3}.
\end{eqnarray*}
Therefore as $m,n\geq n^*(\varepsilon)$,
\begin{eqnarray*}
&&\|{u}_m-{u}_n\|^2_{L^2([0,T];L_\rho^2(\mathbb{R}^d;\mathbb{R}^1))}\nonumber\\
&\leq&\int_0^T\int_{U_{j(\varepsilon)}}|{u}_m(s,x)-{u}_n(s,x)|^2\rho^{-1}(x)dxds+\int_0^T\int_{{U_{j(\varepsilon)}}^c}(2|{u}_m(s,x)|^2+2|{u}_n(s,x)|^2)\rho^{-1}(x)dxds\nonumber\\
&<&\varepsilon.
\end{eqnarray*}
That is to say ${u}_n$ converges strongly in
$L^2([0,T];L_\rho^2(\mathbb{R}^d;\mathbb{R}^1))$. Now using the
equivalence of norm principle, we know as $m$, $n\to\infty$,
\begin{eqnarray}\label{cc}
&&\|Y_s^{t,x,m}-Y_s^{t,x,n}\|^2_{L^2(\Omega\times[t,T];L_\rho^2(\mathbb{R}^d;\mathbb{R}^1))}\nonumber\\
&=&E[\int_t^T\int_{\mathbb{R}^d}|u_m(s,X_s^{t,x})-u_n(s,X_s^{t,x})|^2\rho^{-1}(x)dxds]\nonumber\\
&\leq&C_p\int_t^T\int_{\mathbb{R}^d}|u_m(s,x)-u_n(s,x)|^2\rho^{-1}(x)dxds\longrightarrow0.
\end{eqnarray}
So the claim that $Y_s^{t,x,n}$ converges strongly in
$L^2(\Omega\times[t,T];L_\rho^2(\mathbb{R}^d;\mathbb{R}^1))$
follows. But we know that $Y_s^{t,x}$ is the weak limit of
$Y_s^{t,x,n}$ in
$L^2(\Omega\times[t,T];L_\rho^2(\mathbb{R}^d;\mathbb{R}^1))$,
therefore $Y_s^{t,x,n}$ converges strongly to $Y_s^{t,x}$ in
$L^2(\Omega\times[t,T];L_\rho^2(\mathbb{R}^d;\mathbb{R}^1))$. $\hfill\diamond$\\


Considering the strongly convergent subsequence
$\{Y_\cdot^{t,\cdot,n}\}_{n=1}^\infty$ derived from Lemma \ref{5} and using a standard argument
to BSDE (\ref{c}), we can prove that for arbitrary $m,n$
\begin{eqnarray*}
&&E[\sup_{t\leq s\leq T}\int_{\mathbb{R}^d}|Y_s^{t,x,m}-Y_s^{t,x,n}|^2\rho^{-1}(x)dx]+E[\int_t^T\int_{\mathbb{R}^d}|Z_s^{t,x,m}-Z_s^{t,x,n}|^2\rho^{-1}(x)dxds]\nonumber\\
&\leq&C_pE[\int_t^T\int_{\mathbb{R}^d}|Y_s^{t,x,m}-Y_s^{t,x,n}|^2\rho^{-1}(x)dxds]+C_p\sqrt{E[\int_t^T\int_{\mathbb{R}^d}|Y_s^{t,x,m}-Y_s^{t,x,n}|^2\rho^{-1}(x)dxds]}\nonumber\\
&&\times\sqrt{E[\int_t^T\int_{\mathbb{R}^d}(|f_0(s,x)|^2+|Y_s^{t,x,n}|^{2p}+|Z_s^{t,x,n}|^2)\rho^{-1}(x)dxds]}.
\end{eqnarray*}
So by Condition (H.2) and Lemma \ref{1}, we can conclude that this subsequence $\{Y_\cdot^{t,\cdot,n}\}_{n=1}^\infty$
converges strongly also in $S^{2}([t,T];L_{\rho}^2({\mathbb{R}^{d}};{\mathbb{R}^{1}}))$ and the
corresponding subsequence of $\{Z_\cdot^{t,\cdot,n}\}_{n=1}^\infty$
converges strongly $M^{2}([t,T];L_{\rho}^2({\mathbb{R}^{d}};{\mathbb{R}^{d}}))$ as well.
Certainly the strong convergence limit should be identified with the
weak convergence limit $Z_\cdot^{t,\cdot}$, hence the following
corollary follows without a surprise.
\begin{cor}\label{18}
Let $(Y_\cdot^{t,\cdot},Z_\cdot^{t,\cdot})$ be the solution to BSDE
(\ref{d}) and $(Y_\cdot^{t,\cdot,n},Z_\cdot^{t,\cdot,n})$ be the
subsequence of the solutions to BSDE (\ref{c}), of which
$Y_\cdot^{t,\cdot,n}$ converges strongly to $Y_\cdot^{t,\cdot}$ in
$L^2(\Omega\times[t,T];L_\rho^2(\mathbb{R}^d;\mathbb{R}^1))$, then $(Y_\cdot^{t,\cdot,n},Z_\cdot^{t,\cdot,n})$ also converges strongly to $(Y_\cdot^{t,\cdot},Z_\cdot^{t,\cdot})$
in $S^{2}([t,T];L_{\rho}^2({\mathbb{R}^{d}};{\mathbb{R}^{1}}))\times M^{2}([t,T];L_{\rho}^2({\mathbb{R}^{d}};{\mathbb{R}^{d}}))$.
\end{cor}

As for $Y_s^{t,x}$, we further have
\begin{lem}\label{2}
Under the conditions of Theorem \ref{21},
$E[\int_t^T\int_{\mathbb{R}^d}|Y_s^{t,x}|^{2p}\rho^{-1}(x)dxds]<\infty$ and $Y_s^{t,x}=Y_s^{s,X_s^{t,x}}$ for a.a. $s\in[t,T]$, a.a. $x\in\mathbb{R}^d$ a.s.
\end{lem}
{\em Proof}. First by Lemma \ref{qi045} and Corollary \ref{18}, we have
\begin{eqnarray*}
&&E[\int_t^T\int_{\mathbb{R}^d}|Y_s^{t,x}-Y_s^{s,X_s^{t,x}}|^2\rho^{-1}(x)dxds]\nonumber\\
&\leq&\lim_{n\to\infty}2E[\int_t^T\int_{\mathbb{R}^d}|Y_s^{t,x,n}-Y_s^{t,x}|^2\rho^{-1}(x)dxds]\nonumber\\
&&+\lim_{n\to\infty}2E[\int_t^T\int_{\mathbb{R}^d}|Y_s^{s,X_s^{t,x},n}-Y_s^{s,X_s^{t,x}}|^2\rho^{-1}(x)dxds]\nonumber\\
&\leq&\lim_{n\to\infty}2E[\int_t^T\int_{\mathbb{R}^d}|Y_s^{t,x,n}-Y_s^{t,x}|^2\rho^{-1}(x)dxds]\nonumber\\
&&+\lim_{n\to\infty}C_pE[\sup_{s\leq r\leq T}\int_{\mathbb{R}^d}|Y_r^{s,x,n}-Y_r^{s,x}|^2\rho^{-1}(x)dx]=0.
\end{eqnarray*}
Hence,
\begin{eqnarray}\label{ap}
Y_s^{t,x}=Y_s^{s,X_s^{t,x}}\ {\rm for}\ {\rm a.a.}\ s\in[t,T],\ {\rm a.a.}\
x\in\mathbb{R}^d\ {\rm a.s.}
\end{eqnarray}
If we define $Y_s^{s,x}=u(s,x)$, then by (\ref{ap}) and Lemma \ref{qi045} again we also have
\begin{eqnarray}\label{s}
\lim_{n\to\infty}\int_0^T\int_{\mathbb{R}^d}|u_n(s,x)-u(s,x)|^2\rho^{-1}(x)dxds=0,
\end{eqnarray}
and
\begin{eqnarray*}\label{q}
E[\int_t^T\int_{\mathbb{R}^d}|Y_s^{t,x}-u(s,X_s^{t,x})|^2\rho^{-1}(x)dxds]=0.
\end{eqnarray*}
Therefore, we claim that the strong limit of $u_n(s,x)$ in
$L^2([0,T];L_\rho^2(\mathbb{R}^d;\mathbb{R}^1))$ is $u(s,x)$ and
$Y_s^{t,x}=u(s,X_s^{t,x})$ for a.a. $s\in[t,T]$, $x\in\mathbb{R}^d$
a.s.

By the equivalence of norm principle, to get
$E[\int_t^T\int_{\mathbb{R}^d}|Y_s^{t,x}|^{2p}\rho^{-1}(x)dxds]<\infty$, we only need to
prove
$\int_0^T\int_{\mathbb{R}^d}|u(s,x)|^{2p}\rho^{-1}(x)dxds<\infty$. For this, we first derive from $\lim_{n\to\infty}\int_{0}^{T}\int_{\mathbb{R}^{d}}|u_n(s,x)-u(s,x)|^2\rho^{-1}(x)dxds=0$ a subsequence of
$\{u_n(s,x)\}_{n=1}^\infty$, still denoted by $\{u_n(s,x)\}_{n=1}^\infty$, s.t.
\begin{eqnarray}\label{ao}
u_n(s,x)\longrightarrow u(s,x)\  {\rm and}\
\sup_n|u_n(s,x)|^{2p}<\infty\ \ {\rm for}\ {\rm a.a.}\ s\in[t,T],\
x\in\mathbb{R}^{d}.
\end{eqnarray}
By a similar argument
as in Lemma \ref{10}, for this subsequence $u_n$, we can prove that
for any $\delta>0$,
\begin{eqnarray*}
\lim_{N\to\infty}\sup_n\int_0^T\int_{\mathbb{R}^d}|u_n(s,x)|^{2p-\delta}I_{\{|u_n(s,x)|^{2p-\delta}>N\}}(s,x)\rho^{-1}(x)dxds=0.
\end{eqnarray*}
That is to say that $|u_n(s,x)|^{2p-\delta}$ is uniformly integrable.
Together with $u_n(s,x)\longrightarrow u(s,x)$ for a.a. $s\in[0,T]$,
$x\in\mathbb{R}^{d}$, we have
\begin{eqnarray*}
&&\int_0^T\int_{\mathbb{R}^d}|u(s,x)|^{2p-\delta}\rho^{-1}(x)dxds=\lim_{n\to\infty}\int_0^T\int_{\mathbb{R}^d}|u_n(s,x)|^{2p-\delta}\rho^{-1}(x)dxds\\
&\leq&\sup_{n}\int_0^T\int_{\mathbb{R}^d}|u_n(s,x)|^{2p-\delta}\rho^{-1}(x)dxds\leq
C_p\big(\sup_{n}\int_0^T\int_{\mathbb{R}^d}|u_n(s,x)|^{2p}\rho^{-1}(x)dxds\big)^{{2p-\delta}\over{2p}}\leq
C_p,
\end{eqnarray*}
where the last $C_p<\infty$ is a constant independent of $n$ and
$\delta$.
Then using Fatou lemma to take the limit as $\delta\to0$ in the above inequality, we can get $\int_0^T\int_{\mathbb{R}^d}|u(s,x)|^{2p}\rho^{-1}(x)dxds<\infty$. $\hfill\diamond$\\

Indeed, with Corollary
\ref{18} and Lemma \ref{2}, doing It$\hat {\rm o}$'s formula to
$\psi_M(Y_r^{t,x,n}-Y_r^{t,x})$ and ${\rm
e}^{Kr}\varphi_{n,m}\big(\psi_M(Y_{r}^{t,x})\big)$,  we can further prove that 
$Y_\cdot^{t,\cdot}\in
S^{2p}([t,T];L_{\rho}^{2p}({\mathbb{R}^{d}};{\mathbb{R}^{1}}))$
(To see similar calculations, one can refer to the
argument in the proof of Lemma 3.3 in \cite{zh-zh1}).
\begin{prop}\label{19}
For $(Y_\cdot^{t,\cdot},Z_\cdot^{t,\cdot})$ and
$(Y_\cdot^{t,\cdot,n},Z_\cdot^{t,\cdot,n})$ given in Corollary
\ref{18}, 
$Y_\cdot^{t,\cdot}\in
S^{2p}([t,T];L_{\rho}^{2p}\\({\mathbb{R}^{d}};{\mathbb{R}^{1}}))$.
\end{prop}
Now we are ready to prove the identification of the limiting BSDEs.
\begin{lem}\label{17}
The random field $U$, $Y$ and $Z$ have the following relation:
\begin{eqnarray}\label{ak}
U_s^{t,x}=f(s,X_s^{t,x},Y_s^{t,x},Z_s^{t,x})\ \ {\rm for}\ {\rm
a.a.}\ {\rm s\in[t,T]},\ {\rm x\in\mathbb{R}^d}\ {\rm a.s.}
\end{eqnarray}
\end{lem}
{\em Proof}. 
Let $\mathcal{K}$ be a set in
$\Omega\times[t,T]\times\mathbb{R}^d$ s.t.
$\sup_n|Y_s^{t,x,n}|+\sup_n|Z_s^{t,x,n}|+|f_0(s,X_s^{t,x})|<K$.
Similar to (\ref{ao}), we can find a subsequence of
$\{(Y_s^{t,x,n},Z_s^{t,x,n})\}_{n=1}^\infty$, still denoted by
$\{(Y_s^{t,x,n},Z_s^{t,x,n})\}_{n=1}^\infty$, satisfying
$(Y_s^{t,x,n},Z_s^{t,x,n})\longrightarrow (Y_s^{t,x},Z_s^{t,x})$ and
$\sup_n|Y_s^{t,x,n}|+\sup_n|Z_s^{t,x,n}|<\infty$ for
a.a. $s\in[t,T]$, $x\in\mathbb{R}^{d}$ a.s.
Then it turns out that as $K\to\infty$,
$\mathcal{K}\uparrow\Omega\times[t,T]\times\mathbb{R}^d$. 
Moreover it is easy to see that along the subsequence,
\begin{eqnarray*}
&&E[\int_t^T\int_{\mathbb{R}^d}2(\sup_n|f_n(s,X_s^{t,x},Y_s^{t,x,n},Z_s^{t,x,n})|^2+|f(s,X_s^{t,x},Y_s^{t,x},Z_s^{t,x})|^2)I_\mathcal{K}(s,x)\rho^{-1}(x)dxds]\nonumber\\
&\leq&6C^2E[\int_t^T\int_{\mathbb{R}^d}(|f_0(s,X_s^{t,x})|^2+\sup_n|Y_s^{t,x,n}|^{2p}+\sup_n|Z_s^{t,x,n}|^2)I_\mathcal{K}(s,x)\rho^{-1}(x)dxds]\nonumber\\
&&+6C^2E[\int_t^T\int_{\mathbb{R}^d}(|f_0(s,X_s^{t,x})|^2+|Y_s^{t,x}|^{2p}+|Z_s^{t,x}|^2)I_\mathcal{K}(s,x)\rho^{-1}(x)dxds]<\infty.
\end{eqnarray*}
Thus, we can apply Lebesgue's dominated convergence theorem to the
following calculation:
\begin{eqnarray}\label{zz17}
&&\lim_{n\to\infty}E[\int_t^T\int_{\mathbb{R}^d}|f_n(s,X_s^{t,x},Y_s^{t,x,n},Z_s^{t,x,n})I_\mathcal{K}(s,x)-f(s,X_s^{t,x},Y_s^{t,x},Z_s^{t,x})I_\mathcal{K}(s,x)|^2\rho^{-1}(x)dxds]\nonumber\\
&=&E[\int_t^T\int_{\mathbb{R}^d}\lim_{n\to\infty}|f_n(s,X_s^{t,x},Y_s^{t,x,n},Z_s^{t,x,n})-f(s,X_s^{t,x},Y_s^{t,x},Z_s^{t,x})|^2I_\mathcal{K}(s,x)\rho^{-1}(x)dxds]\nonumber\\
&\leq&2E[\int_t^T\int_{\mathbb{R}^d}\lim_{n\to\infty}|f_n(s,X_s^{t,x},Y_s^{t,x,n},Z_s^{t,x,n})-f(s,X_s^{t,x},Y_s^{t,x,n},Z_s^{t,x,n})|^2I_\mathcal{K}(s,x)\rho^{-1}(x)dxds]\nonumber\\
&&+2E[\int_t^T\int_{\mathbb{R}^d}\lim_{n\to\infty}|f(s,X_s^{t,x},Y_s^{t,x,n},Z_s^{t,x,n})-f(s,X_s^{t,x},Y_s^{t,x},Z_s^{t,x})|^2I_\mathcal{K}(s,x)\rho^{-1}(x)dxds].\nonumber\\
\end{eqnarray}
Since $Y_s^{t,x,n}\longrightarrow Y_s^{t,x}$ for a.a. $s\in[t,T]$,
$x\in\mathbb{R}^{d}$ a.s., there exists a $N(s,x,\omega)$ s.t. when
$n\geq N(s,x,\omega)$, $|Y_s^{t,x,n}|\leq|Y_s^{t,x}|+1$. So taking
$n\geq\max\{N(s,x,\omega),\ |Y_s^{t,x}|+1\}$, we have
$f_n(s,X_s^{t,x},Y_s^{t,x,n},Z_s^{t,x,n})=f(s,X_s^{t,x},{\inf(n,|Y_s^{t,x,n}|)\over|Y_s^{t,x,n}|}Y_s^{t,x,n},Z_s^{t,x,n})=f(s,X_s^{t,x},Y_s^{t,x,n},Z_s^{t,x,n})$.
That is to say\\
$\lim_{n\to\infty}|f_n(s,X_s^{t,x},Y_s^{t,x,n},Z_s^{t,x,n})-f(s,X_s^{t,x},Y_s^{t,x,n},Z_s^{t,x,n})|^2=0$
for a.a. $s\in[t,T]$, $x\in\mathbb{R}^{d}$ a.s. On the other hand,
$\lim_{n\to\infty}|f(s,X_s^{t,x},Y_s^{t,x,n},Z_s^{t,x,n})-f(s,X_s^{t,x},Y_s^{t,x},Z_s^{t,x})|^2=0$
for a.a. $s\in[t,T]$, $x\in\mathbb{R}^{d}$ a.s. is obvious due to
the continuity of $(y,z)\to f(s,x,y,z)$.

Therefore by (\ref{zz17}),
$f_n(s,X_s^{t,x},Y_s^{t,x,n},Z_s^{t,x,n})I_\mathcal{K}(s,x)=U_s^{t,x,n}I_\mathcal{K}(s,x)$
converges strongly to\\
$f(s,X_s^{t,x},Y_s^{t,x},Z_s^{t,x})I_\mathcal{K}(s,x)$ in
$L^2_\rho(\Omega\times[t,T]\times\mathbb{R}^d;\mathbb{R}^1)$, but
$U_s^{t,x,n}I_\mathcal{K}(s,x)$ converges weakly to
$U_s^{t,x}I_\mathcal{K}(s,x)$ in
$L^2_\rho(\Omega\times[t,T]\times\mathbb{R}^d;\mathbb{R}^1)$, so
$f(s,X_s^{t,x},Y_s^{t,x},Z_s^{t,x})I_\mathcal{K}(s,x)=U_s^{t,x}I_\mathcal{K}(s,x)$
for a.a. $r\in[t,T]$, $x\in\mathbb{R}^{d}$ a.s. The lemma follows
when $K\to\infty$. $\hfill\diamond$\\

{\em Proof of Theorem \ref{21}}. With Proposition
\ref{19} and Lemma \ref{17}, the existence of solutions to BSDE (\ref{t}) is easy to
see. Now we prove the uniqueness. If there is another solution
$(\tilde{Y}_s^{t,x},\tilde{Z}_s^{t,x})$ to BSDE (\ref{t}), then for
a.a. $x\in\mathbb{R}^d$,
$
(Y_s^{t,x}-\tilde{Y}_s^{t,x},Z_s^{t,x}-\tilde{Z}_s^{t,x})$ satisfies
\begin{eqnarray}\label{aa}
Y_s^{t,x}-\tilde{Y}_s^{t,x}=\int_s^T\big(f(r,X_r^{t,x},Y_r^{t,x},Z_r^{t,x})-f(r,X_r^{t,x},\tilde{Y}_r^{t,x},\tilde{Z}_r^{t,x})\big)dr-\int_s^T\langle
Z_r^{t,x}-\tilde{Z}_r^{t,x},dW_r\rangle.\nonumber
\end{eqnarray}
Applying It$\hat {\rm o}$'s formula to
$|Y_s^{t,x}-\tilde{Y}_s^{t,x}|^2$, by the stochastic Fubini theorem
and Conditions ${\rm (H.3)}^*$ and (H.4), we have
\begin{eqnarray*}
&&E[\int_{\mathbb{R}^d}|Y_s^{t,x}-\tilde{Y}_s^{t,x}|^2\rho^{-1}(x)dx]+E[\int_{s}^{T}\int_{\mathbb{R}^d}|Z_r^{t,x}-\tilde{Z}_r^{t,x}|^2\rho^{-1}(x)dxdr]\nonumber\\
&\leq&2L^2E[\int_{s}^{T}\int_{\mathbb{R}^d}|Y_r^{t,x}-\tilde{Y}_r^{t,x}|^2\rho^{-1}(x)dxdr]+{1\over2}E[\int_{s}^{T}\int_{\mathbb{R}^d}|Z_r^{t,x}-\tilde{Z}_r^{t,x}|^2\rho^{-1}(x)dxdr].
\end{eqnarray*}
By Gronwall's inequality, the uniqueness of the solution to BSDE (\ref{t}) follows immediately. $\hfill\diamond$\\


By the stochastic flow $X_r^{s,X^{t,x}_s}=X^{t,x}_r$ for $t\leq s\leq r\leq T$ and the uniqueness of solution of BSDE (\ref{t}), following a similar argument as Proposition 3.4 in \cite{zh-zh1} we have
\begin{cor}\label{4}
Under the conditions of Theorem \ref{21}, let $(Y_s^{t,x},Z_s^{t,x})$ be the solution of BSDE
(\ref{t}), then
\begin{eqnarray*}\label{v}
Y_s^{t,x}=Y_s^{s,X_s^{t,x}},\ Z_s^{t,x}=Z_s^{s,X_s^{t,x}}\ \ {\rm
for}\ {\rm any}\ s\in[t,T],\ {\rm a.a.}\ x\in\mathbb{R}^d\ {\rm
a.s.}
\end{eqnarray*}
\end{cor}

\section{The PDEs}
\setcounter{equation}{0}

Now we make use of the results for BSDE (\ref{t}) to give the
probabilistic representation to PDEs with p-growth coefficients.
Actually the solution of BSDE in the $\rho$-weighted $L^2$ space
gives the unique weak solution of its corresponding PDE (\ref{z}).\\

{\em Proof of Theorem \ref{20}}. Using Corollary \ref{18}, we first
prove the relationship between $(Y,Z)$ and $u$, when we take $u(t,x)=Y_t^{t,x}$. Having proved Lemma \ref{2}, we
only need to prove that $(\sigma^*\nabla u)(s,X^{t,x}_s)=Z_s^{t,x}$
for a.a. $s\in[t,T]$, $x\in\mathbb{R}^{d}$ a.s. This can be deduced
from Corollary \ref{4} and the strong convergence of $Z_\cdot^{t,\cdot,n}$
to $Z_\cdot^{t,\cdot}$ in
$L^2(\Omega\times[t,T];L_\rho^2(\mathbb{R}^d;\mathbb{R}^1))$ by the
similar argument as in Proposition 4.2 in \cite{zh-zh1}.

We then prove that $u(t,x)$ defined above is the unique weak solution of PDE
(\ref{z}). We still start from PDE (\ref{f}). Let $u^n(s,x)$ be the
weak solution of PDE (\ref{f}). Then by the definition for the weak
solution of PDE, we know $(u_n,\sigma^*\nabla u_n)\in
L^{2}([0,T];L_{\rho}^2({\mathbb{R}^{d}};{\mathbb{R}^{1}}))\times
L^{2}([0,T];L_{\rho}^2({\mathbb{R}^{d}};{\mathbb{R}^{d}}))$ and for
an arbitrary $\varphi\in C_c^{\infty}(\mathbb{R}^d;\mathbb{R}^1)$,
\begin{eqnarray}\label{qi16}
&&\int_{\mathbb{R}^{d}}u_n(t,x)\varphi(x)dx-\int_{\mathbb{R}^{d}}u_n(T,x)\varphi(x)dx-{1\over2}\int_{t}^{T}\int_{\mathbb{R}^{d}}\big((\sigma^*\nabla u_n)(s,x)\big)^*(\sigma^*\nabla\varphi)(x)dxds\nonumber\\
&&-\int_{t}^{T}\int_{\mathbb{R}^{d}}u_n(s,x)div\big((b-\tilde{A})\varphi\big)(x)dxds\nonumber\\
&=&\int_{t}^{T}\int_{\mathbb{R}^{d}}f_n\big(s,x,u_n(s,x),(\sigma^*\nabla
u_n)(s,x)\big)\varphi(x)dxds.
\end{eqnarray}
We can prove along a subsequence that each term of (\ref{qi16})
converges to the corresponding term of (\ref{ae}). By (\ref{s}), we
know that $u_n$ converges strongly to $u$ in
$L^2_\rho([0,T]\times\mathbb{R}^d;\mathbb{R}^1)$, thus 
$u_n$ also converges weakly. Moreover,
$\sup_{x\in\mathbb{R}^d}(|div\big((b-\tilde{A})\varphi\big)(x)|)<\infty$
and $\rho$ is a continuous function in $\mathbb{R}^d$, so it is
obvious that
\begin{eqnarray*}
\lim_{n\to\infty}\int_{t}^{T}\int_{\mathbb{R}^{d}}u_n(s,x)div\big((b-\tilde{A})\varphi\big)(x)dxds=\int_{t}^{T}\int_{\mathbb{R}^{d}}u(s,x)div\big((b-\tilde{A})\varphi\big)(x)dxds.
\end{eqnarray*}
Also it is easy to see that
\begin{eqnarray*}
&&\lim_{n\to\infty}{1\over2}\int_{t}^{T}\int_{\mathbb{R}^{d}}\big((\sigma^*\nabla u_n)(s,x)\big)^*(\sigma^*\nabla\varphi)(x)dxds\nonumber\\
&=&\lim_{n\to\infty}-{1\over2}\int_{t}^{T}\int_{\mathbb{R}^{d}}u_n(s,x)div(\sigma\sigma^*\nabla\varphi)(x)\rho(x)\rho^{-1}(x)dxds\nonumber\\
&=&-{1\over2}\int_{t}^{T}\int_{\mathbb{R}^{d}}u(s,x)div(\sigma\sigma^*\nabla\varphi\sigma)(x)\rho(x)\rho^{-1}(x)dxds\nonumber\\
&=&{1\over2}\int_{t}^{T}\int_{\mathbb{R}^{d}}\big((\sigma^*\nabla
u)(s,x)\big)^*(\sigma^*\nabla\varphi)(x)dxds.
\end{eqnarray*}
Also we have proved that $f_n(s,X_s^{t,x},Y_s^{t,x,n},Z_s^{t,x,n})$
converges weakly to $f(s,X_s^{t,x},Y_s^{t,x},Z_s^{t,x})$ in
$L_\rho^2(\Omega\times[t,T]\times\mathbb{R}^d;\mathbb{R}^1)$. In
fact we can follow the same procedure as in the proof of Lemma
\ref{17} to prove $f_n\big(s,x,u_n(s,x),(\sigma^*\nabla
u_n)(s,x)\big)$ converges weakly to
$f\big(s,x,u(s,x),(\sigma^*\nabla u)(s,x)\big)$ in
$L_\rho^2([t,T]\times\mathbb{R}^d;\mathbb{R}^1)$. So we have
\begin{eqnarray*}
&&\lim_{n\to\infty}\int_t^T\int_{\mathbb{R}^d}f_n\big(s,x,u_n(s,x),(\sigma^*\nabla
u_n)(s,x)\big)\varphi(x)dxds\\
&=&\int_t^T\int_{\mathbb{R}^d}f\big(s,x,u(s,x),(\sigma^*\nabla
u)(s,x)\big)\varphi(x)dxds.
\end{eqnarray*}
For any $t\in[0,T]$,
$\lim_{n\to\infty}\int_{\mathbb{R}^{d}}u_n(t,x)\varphi(x)dx=\int_{\mathbb{R}^{d}}u(t,x)\varphi(x)dx$
can be proved as follows using Corollary \ref{18}:
\begin{eqnarray*}
\lim_{n\to\infty}|\int_{\mathbb{R}^{d}}(u_n(t,x)-u(t,x))\varphi(x)dx|^2
&\leq&\lim_{n\to\infty}C_pE[\int_{\mathbb{R}^{d}}|u_n(t,X_t^{0,x})-u(t,X_t^{0,x})|^2\rho^{-1}(x)dx]\nonumber\\
&\leq&\lim_{n\to\infty}C_pE[\sup_{0\leq t\leq
T}\int_{\mathbb{R}^{d}}|Y_t^{0,x,n}-Y_t^{0,x}|^2\rho^{-1}(x)dx]=0.\nonumber
\end{eqnarray*}
Here the convergence in the $S^{2p}$ space gives us a strong result about the convergence of $\int_{\mathbb{R}^d}u_n(t,x)\varphi(x)dx\\\longrightarrow\int_{\mathbb{R}^d}u(t,x)\varphi(x)dx$ uniformly in $t$ as $n\to\infty$.
Therefore we can prove (\ref{ae}) is satisfied for all $t\in[0,T]$.
That is to say $u(t,x)$ is a weak solution of PDE (\ref{z}).

The uniqueness of PDE (\ref{z}) can be derived from the uniqueness
of BSDE (\ref{t}). Let $u$ be a solution of PDE (\ref{z}). Define
$F(s,x)=f\big(s,x,u(s,x),(\sigma^*\nabla u)(s,x)\big)$. Since
$u$ is a solution, so
$\int_{0}^{T}\int_{\mathbb{R}^d}\big(|u(s,x)|^{2p}+|(\sigma^*\nabla
u)(s,x)|^2\big)\rho^{-1}(x)dxds<\infty$ and
\begin{eqnarray}\label{zz14}
&&\int_{0}^{T}\int_{\mathbb{R}^d}|F(s,x)|^2\rho^{-1}(x)dxds\nonumber\\
&\leq&C_p\int_{0}^{T}\int_{\mathbb{R}^d}\big(|f_0(s,x)|^2+|u(s,x)|^{2p}+|(\sigma^*\nabla
u)(s,x)|^2\big)\rho^{-1}(x)dxds<\infty.
\end{eqnarray}
Then we get a PDE with the generator $F\in L^{2}([0,T];L_{\rho}^2({\mathbb{R}^{d}};{\mathbb{R}^{1}}))$. For this generator $F$, we claim that $(Y_s^{t,x},Z_s^{t,x})\triangleq (u(s,X_s^{t,x}),(\sigma^*\nabla u)(s,X^{t,x}_s))$ solves the following linear BSDE for a.a. $x\in\mathbb{R}^d$ with probability one:
\begin{eqnarray}\label{zz15}
Y_s^{t,x}=h(X_{T}^{t,x})+\int_{s}^{T}F(r,X_{r}^{t,x})dr-\int_{s}^{T}\langle
Z^{t,x}_r,dW_r\rangle.
\end{eqnarray}
First we use the mollifier to smootherize $h$ and $F$, then we get two smootherized sequences $h^m$ and $F^m$ such that $h^m(\cdot)\longrightarrow h(\cdot)$ and
$F^m(s,\cdot)\longrightarrow F(s,\cdot)$ in
$L_{\rho}^2({\mathbb{R}^{d}};{\mathbb{R}^{1}})$ respectively. Denote
by ${u}_m(t,x)$ the solution
of PDE on $[0,T]$ with terminal value $h^m(x)$ and generator $F^m(s,x)$ and by $(Y_{s,m}^{t,x},Z_{s,m}^{t,x})$ the solution of BSDE with terminal value $h^m(X_T^{t,x})$ and generator $F(s,X_s^{t,x})$, then following classical results of Pardoux and Peng \cite{pa-pe2}, we
have ${Z}_{t,m}^{t,x}=\sigma^*\nabla{u}_m(t,x)$, and ${Y}_{s,m}^{t,x}={u}_m(s,X_s^{t,x})={Y}_{s,m}^{s,X_s^{t,x}}$,
${Z}_{s,m}^{t,x}=\sigma^*\nabla{u}_m(s,X_s^{t,x})={Z}_{s,m}^{s,X_s^{t,x}}$.
But by standard estimates $({Y}_{s,m}^{t,x},{Z}_{s,m}^{t,x})$ is a Cauchy sequence in $L^2(\Omega\times[t,T];L_{\rho}^2({\mathbb{R}^{d}};{\mathbb{R}^{1}}))\times
L^2(\Omega\times[t,T];L_{\rho}^2({\mathbb{R}^{d}};{\mathbb{R}^{d}}))$. By equivalence of norm principle, $u_m(s,x)$ is also a Cauchy sequence in $\mathcal{H}$, where $\mathcal{H}$ is the set of random fields $\{w(s,x);\
s\in[0,T],\ x\in\mathbb{R}^{d}\}$ such that $(w,\sigma^*\nabla w)\in
L^2([0,T];L_{\rho}^2({\mathbb{R}^{d}};{\mathbb{R}^{1}}))\times
L^2([0,T];L_{\rho}^2({\mathbb{R}^{d}};{\mathbb{R}^{d}}))$ with
the norm $\sqrt{E[\int_{0}^{T}\int_{\mathbb{R}^{d}}(|w(s,x)|^2
+|(\sigma^*\nabla)w(s,x)|^2)\rho^{-1}(x)dxds)}<\infty$. So there exists ${u}\in\mathcal{H}$ such that
$({u}_m,\sigma^*\nabla
{u}_m)\rightarrow({u},\sigma^*\nabla{u})$
in
$L^2([0,T];L_{\rho}^2({\mathbb{R}^{d}};{\mathbb{R}^{1}}))\times
L^2([0,T];L_{\rho}^2({\mathbb{R}^{d}};{\mathbb{R}^{d}}))$ due to the completeness of $\mathcal{H}$. By the equivalence of norm principle again, we know that $(Y_s^{t,x},Z_s^{t,x})\triangleq (u(s,X_s^{t,x}),(\sigma^*\nabla u)(s,X^{t,x}_s))$ is the limit of Cauchy sequence of $({Y}_{s,m}^{t,x},{Z}_{s,m}^{t,x})$. Now it is
easy to pass the limit as $m\rightarrow\infty$ on the BSDE which $({Y}_{s,m}^{t,x},{Z}_{s,m}^{t,x})$ satisfies
and conclude that $(Y_s^{t,x},Z_s^{t,x})$ is a solution
of BSDE (\ref{zz15}).

Noting the
definition of $F(s,x)$, $Y_s^{t,x}$ and $Z_s^{t,x}$, we have that
$(Y_s^{t,x},Z_s^{t,x})$ solves
BSDE (\ref{t}) for a.a. $x\in\mathbb{R}^d$ with probability one.
Moreover, by Lemma
\ref{qi045},
\begin{eqnarray*}
&&E[\int_{t}^{T}\int_{\mathbb{R}^d}(|Y_s^{t,x}|^{2p}+|Z_s^{t,x}|^2)\rho^{-1}(x)dxds]\nonumber\\
&\leq&C_p\int_{t}^{T}\int_{\mathbb{R}^d}|u(s,x)|^{2p}+|(\sigma^*\nabla
u)(s,x)|^2\rho^{-1}(x)dxds<\infty.
\end{eqnarray*}
As Proposition \ref{19}, we can further deduce
that $Y_\cdot^{t,\cdot}\in
S^{2p}([t,T];L_{\rho}^{2p}({\mathbb{R}^{d}};{\mathbb{R}^{1}}))$ and
therefore $(Y_s^{t,x},Z_s^{t,x})$ is a solution of BSDE (\ref{t}).
If there is another solution $\hat{u}$ to PDE (\ref{z}), then by the
same procedure, we can find another solution
$(\hat{Y}_s^{t,x},\hat{Z}_s^{t,x})$ to BSDE (\ref{t}), where
\begin{eqnarray*}
\hat{Y}_s^{t,x}=\hat{u}(s,X_s^{t,x})\ {\rm and}\
\hat{Z}_s^{t,x}=(\sigma^*\nabla \hat{u})(s,X^{t,x}_s).
\end{eqnarray*}
By Theorem \ref{21}, the solution of BSDE (\ref{t}) is unique.
Therefore
\begin{eqnarray*}
Y_s^{t,x}=\hat{Y}_s^{t,x}\ {\rm for}\ {\rm a.a.}\ s\in[t,T],\ {\rm a.a.}\ x\in\mathbb{R}^{d}\ {\rm a.s.}
\end{eqnarray*}
In particular, when $t=0$,
\begin{eqnarray*}
Y_s^{0,x}=\hat{Y}_s^{0,x}\ {\rm for}\ {\rm a.a.}\ s\in[0,T],\
x\in\mathbb{R}^{d}\ {\rm a.s.}
\end{eqnarray*}
By Lemma \ref{qi045} again,
\begin{eqnarray*}
\int_{0}^{T}\int_{\mathbb{R}^d}|u(s,x)-\hat{u}(s,x)|^2\rho^{-1}(x)dxds
\leq
C_pE[\int_{0}^{T}\int_{\mathbb{R}^d}|Y_s^{0,x}-\hat{Y}_s^{0,x}|^2)\rho^{-1}(x)dxds]
=0.
\end{eqnarray*}
So $u(s,x)=\hat{u}(s,x)$ for a.a. $s\in[0,T]$, $x\in\mathbb{R}^{d}$
a.s. The uniqueness is proved. The uniqueness implies that for any selection $u$ in the equivalence class of solution of the PDE (\ref{z}), $u(s,x)=Y_s^{s,x}$ for a.a. $s\in[0,T]$, $x\in\mathbb{R}^d$. Moreover, noting that $(u(s,X_s^{t,x}),\sigma^*\nabla u(s,X_s^{t,x}))$ solves the BSDE (\ref{t}) and using the uniqueness of solution of BSDE (\ref{t}) in the equivalence class, we have (\ref{ce}) for any representative $Y$ in the equivalence class of the solution of BSDE (\ref{t}).
$\hfill\diamond$\\

{\bf Acknowledgements}. We would like to acknowledge useful
conversations with K.D. Elworthy, C.R. Feng, J. Lorinczi, K.N. Lu,
S.G. Peng and S.J. Tang. QZ would like to thank the Department of
Mathematical Sciences of Loughborough University for appointing him
as a Research Associate. He wishes to acknowledge their financial
support to the project through the appointment and partial financial
support of the National Basic Research Program of China (973
Program) with Grant No.2007CB814904. HZ would like to thank K.N. Lu
for inviting him to visit Brigham-Young University under their
programme of Special Year on Stochastic Dynamics and S.J. Tang for
inviting him to visit Laboratory of Mathematics for Nonlinear
Science, Fudan University. We are grateful to the referee for his/her constructive comments.

\end{document}